\documentclass[11pt,intlimits]{article}

\usepackage{amsbsy,amsmath,amssymb,amsfonts}
\usepackage{epsfig}

\newtheorem{thm}{Theorem}[section]

\newtheorem{defi}{Definition}[section]
\newtheorem{lemma}{Lemma}[section]
\newtheorem{rem}{Remark}[section]
\newtheorem{cond}{Condition}[section]

\newcommand{\paragrafo}[1]{\section{#1}\setcounter{equation}{0}}
\newcommand{\subparagrafo}[1]{\subsection{#1}}

\def\R{\mathbb{R}}
\def\Z{\mathbb{Z}}
\def\de{\partial}

\def\G{\Gamma}

\def\e{\varepsilon}
\def\th{\vartheta}
\def\Th{\Theta}
\def\z{\zeta}
\def\h{\eta}
\def\k{\kappa}
\def\l{\lambda}

\def\m{\mu}

\def\r{\varrho}
\def\s{\sigma}

\def\o{\omega}

\def\K{{\mathcal K}}
\def\M{{\mathcal M}}

\def\Ll{{\mathcal L}}
\def\C{{\mathcal C}}
\def\B{{\mathcal B}}

\def\Ss{{\mathcal S}}
\def\Ff{{\mathcal F}}

\def\Pol{{\mathcal P}}

\def\T{{\mathcal T}}

\def\ee{{\rm e}}

\def\leq{\leqslant}
\def\geq{\geqslant}
\def\bg{{\boldsymbol{\gamma}}}
\def\ba{{\boldsymbol{\alpha}}}
\def\bb{{\boldsymbol{\beta}}}
\def\bmu{{\boldsymbol{\mu}}}

\def\bm{{\mathbf m}}
\def\bj{{\mathbf j}}

\def\bk{{\mathbf k}}
\def\bx{{\mathbf x}}
\def\bX{{\mathbf X}}
\def\by{{\mathbf y}}
\def\bz{{\mathbf z}}

\def\cc{\mathbf c}

\def\bxi{{\mathbf \xi}}

\def\btt{{\mathbf t}}

\def\dd{\sqrt{D}}
\def\hj{\eta_\bj}

\def\asb{(\bx)}

\def\uxj{u(\bx_\bj)}
\def\st{\mathrm {st}\,}
\def\ST{\mathrm {ST}\,}
\def\supp{{\mathrm {supp}}\,}
\def\ds{\displaystyle}
\def\ov{\overline}

\title{Approximate Approximations from scattered data}

\author{{\small F. LANZARA$^{\mbox{\tiny 1}}$ , V. MAZ'YA$^{\mbox{\tiny 2}}$ ,
G. SCHMIDT$^{\mbox{\tiny 3}}$}
}
\date{}
\begin{document}

\maketitle

\baselineskip=0.9\normalbaselineskip
\vspace{-3pt}

\hspace*{1cm}
\parbox{10cm}{\begin{flushleft}
{\footnotesize\em
\begin{itemize}
\item[$^{\mbox{\tiny\rm 1}}$]Dipartimento di Matematica, Universit\`a
``La Sapienza'',\\
Piazzale Aldo Moro 2, 00185 Roma, Italy\\
\texttt{\rm lanzara\symbol{'100}mat.uniroma1.it}
\item[$^{\mbox{\tiny\rm 2}}$]Department of Mathematics, University of
Link\"oping, \\ 581 83 Link\"oping, Sweden\\
\texttt{\rm vlmaz\symbol{'100}mai.liu.se }
\item[$^{\mbox{\tiny\rm 3}}$]Weierstrass Institute for Applied Analysis and
Stochastics, \\  Mohrenstr. 39,
10117 Berlin, Germany \\
\texttt{\rm schmidt\symbol{'100}wias-berlin.de}
\end{itemize}
}
\end{flushleft}}

\begin{abstract}
The aim of this paper is to extend the 
approximate quasi-interpolation on a uniform grid
by dilated shifts of a smooth and rapidly decaying function
on a uniform grid to scattered data quasi-interpolation. 
It is shown that high order approximation of smooth functions 
up to some prescribed accuracy  is possible, if the basis functions,
which are centered at the scattered nodes, are multiplied
by suitable polynomials such that
their sum is an approximate partition of unity.   
For Gaussian functions we propose a method to construct
the  approximate partition of unity and describe the application 
of the new quasi-interpolation approach
to
the cubature of multi-dimensional integral operators.
\end{abstract}

\paragrafo{Introduction}\label{S:intro}
\setcounter{equation}{0}
\vspace{.2in}
The approximation of multivariate functions  from scattered data 
is an important theme in numerical mathematics. 
One of the methods to attack this problem is quasi-interpolation.
One takes values $u(\bx_\bj)$ of a
function $u$ on a set of nodes $\{\bx_\bj\}$ 
and constructs an approximant of $u$ by linear
combinations
\[
\sum u(\bx_\bj) \h_\bj(\bx) \, ,
\]   
where $\h_\bj(\bx)$ is a set of basis
functions. 
Using quasi-interpolation there is no need to solve large algebraic systems.
The approximation properties of quasi-interpolants
in the case that $\bx_\bj$ are the nodes of a uniform grid
are well-understood.
For example, the quasi-interpolant
\begin{equation} \label{d:simple}
\sum_{\bj \in \Z^n} u(h \bj) \varphi \Big(\frac{\bx-h \bj }{h} \Big) 
\end{equation}
can be studied via the theory of principal shift-invariant
spaces, which has been developed  in several articles by de Boor, 
DeVore and Ron
(see {\it  e.g.} \cite{Bo}, \cite{BR}).
Here $\varphi$ is supposed to be a compactly supported 
or rapidly decaying function.
Based on the Strang-Fix condition for $\varphi$, 
which is equivalent to polynomial reproduction,
convergence  and approximation orders for several classes of
basis functions were obtained (see also Schaback/Wu \cite{SWu},
Jetter/Zhou \cite{JZ}).
Scattered data quasi-interpolation by functions, which reproduce
polynomials,  
has been studied
by Buhmann, Dyn, Levin in \cite{BDL} and Dyn, Ron in \cite{DR}
(see also \cite{Wu} for further references).

In order to extend the
quasi-interpolation \eqref{d:simple}
to general classes of approximating functions, 
another concept of approximation procedures, called  {\it Approximate Approximations},
was proposed in \cite{Ma1} and \cite{Ma2}.
These procedures have the common feature, that
they are accurate without
being convergent in a rigorous sense.
Consider, for example,
the quasi-interpolant on the uniform grid
\begin{equation}\label{1.1}
\M\,u(\bx)=D^{-n/2} \sum_{\bj \in \Z^n} u(h\bj)\,
\h{\Big(\frac{\bx -h \bj}{h \sqrt{D}}\Big)} \, ,
\end{equation}
where $\h$ is sufficiently smooth and of rapid decay,
$h$ and $D$ are two positive parameters.
It was shown that if $\Ff \eta -1$ has a zero of order
 $N$ at the origin ($\Ff \eta$ denotes the Fourier transform of $\h$),
then $\M\,u$ approximates $u$ pointwise 
\begin{equation}\label{est1bis}
|\M u(\bx)-u(\bx)|\leq c_{N,\h}\,(h \dd)^{N} \sup_{\R^n}|\nabla_{N}u|+
\e \, |\nabla_{N-1}u(\bx)| 
\end{equation}
with a constant $c_{N,\h}$ not depending on $u$, $h$, and $D$, and 
$\e$ can be made arbitrarily small if $D$ is sufficiently large
(see \cite{MS2}, \cite{MS3}).
In general, there is no convergence of the {\it approximate quasi-interpolant}
$\M u(\bx)$ to $u(\bx)$ as $h \to 0$.
However, one can fix 
$D$ such that up to any prescribed accuracy
$\M u$ approximates $u$ with order $O(h^N)$.
The lack of convergence as $h \to 0$, which is  even not
perceptible in numerical computations for appropriately chosen $D$, 
is compensated by a greater
flexibility in the choice of approximating functions $\h$.
In applications, this flexibility enables one
to obtain simple and accurate 
formulae for values of various integral and
pseudo-differential operators of mathematical physics
(see  \cite{MS1}, \cite{MS5}, \cite{MSW} and the review
paper \cite{Sc})
and to develop explicit
semi-analytic time-marching algorithms for initial
boundary value problems for linear and non linear evolution equations
(\cite{KM1}, \cite{KM2}).

The approximate quasi-interpolation approach 
was extended to nonuniform grids up to now in
two directions.
The case that the set of nodes is
a smooth image of a  uniform grid have been studied  in \cite{MS4}. 
It was shown   
that formulae similar to (\ref{1.1}) 
preserve the basic properties of approximate quasi-interpolation.
A similar result for quasi-interpolation on piecewise uniform grids
was obtained in \cite {IMS}.

It is the purpose of 
the present paper to generalize the method of
approximate quasi-interpolation
to functions with values given on
a rather general  grid.
We start with a simple quasi-interpolant
for a set of nodes close to a uniform grid of size $h$ in the
sense, that for some positive constant $\kappa$ and  any $\bj \in \Z^n$ there exists
at least one node $\bx_\bj$ with $|\bx_\bj - h \bj| < \kappa h$.
Then
under some additional assumption on the nodes we construct
a quasi-interpolant with gridded centers
\[
{\mathbb M} u(\bx)= D^{-n/2} \sum_{\bj \in \Z^n} \Lambda_\bj(u)  
\h\Big( \frac{\bx-h\, \bj}{h \dd}\Big) \, .
\]
Here $\Lambda_\bj$ are linear functionals  of the data at
a finite number of nodes around  $\bx_\bj$.
It can be shown that estimate \eqref{est1bis} remains true
for ${\mathbb M} u$ 
under the same assumptions on the
function $\h$.

In order to treat more general distributions of the  nodes $\bx_\bj$
we modify the approximating functions.
More precisely, we 
consider approximations of the form
\begin{equation}\label{e:nonuniform}
Mu(\bx)=\sum_\bj \uxj \Pol_\bj(\bx) 
\h{\Big(\frac{\bx -\bx_\bj}{h_\bj}\Big)}
\end{equation}
with some polynomials $\Pol_\bj$.
We show that  one can achieve the approximation of
$u$ with arbitrary order $N$
up to a small saturation error,
as long as an "approximate partition of unity" 
$\ds \Big\{\widetilde \Pol_\bj(\bx)\}\h\Big(\frac{\bx -\bx_\bj}{h_\bj}\Big) \Big\}$
with other polynomials $\widetilde \Pol_\bj$ exists.
Here we mean that for any $\e >0$ one can find polynomials such that 
\begin{equation*}
\Big|\sum_\bj \widetilde \Pol_\bj(\bx) 
\h{\Big(\frac{\bx -\bx_\bj}{h_\bj}\Big)} - 1 \Big| < \e \, .
\end{equation*}
Then one can choose the polynomials $\Pol_\bj$ in
\eqref{e:nonuniform} such that 
\[
|M u(\bx)-u(\bx)|\leq C \, \sup_{\bj} h^N_\bj \, \|\nabla_{N}u\|_{L_\infty} +
\e \, |u(\bx)| \, .
\]
This estimate is valid as long as 
\[
\sum_\bj \h\Big(\frac{\bx -\bx_\bj}{h_\bj}\Big) \ge c >0
\]
and $\h$ is sufficiently smooth and of rapid decay, but is not
subjected to additional requirements as the Strang-Fix condition.
Moreover, we propose a method  to construct the polynomials 
such that the series
\[
\sum_\bj \widetilde \Pol_\bj(\bx) \,\ee^{- |\bx -\bx_j|^2/h^2_\bj D}
\]
approximates the constant 
function  $1$  up to an arbitrary prescribed accuracy.
This method does not require
solving a large system of linear equations. Instead, in order to obtain
the local representation of the partition of unity, one has to solve 
a small number of approximation problems, which are reduced
to linear systems of moderate size. 

By a suitable choice
of $\h$  it is possible
to obtain explicit semi-analytic or other efficient
approximation
formulae for multi-dimensional integral and pseudo-differential
operators which are based on the quasi-interpolant \eqref{e:nonuniform}.
So the cubature of those integrals, which is one of the applications
of the approximate quasi-interpolation on uniform grids, can be
carried over to the case when the
integral operators are applied to functions 
given at scattered nodes.

We give a simple example of
formula 
(\ref{e:nonuniform}).
Let $\{x_i\}$ be a sequence of points on $\R$ such that $0<x_{i+1}-x_i\leq 1$.
Consider a sequence  of functions $\z_j$ on $\R$ supported by a fixed
neighborhood of the origin.
Suppose that the sequence $\{\z_j(x-x_j)\}$
forms an approximate
partition of unity on $\R$,
\[
|1- \sum_j \z_j(x-x_j)|<\e\>\> .
\]
One can easily see that the quasi-interpolant
\begin{equation*}
\begin{split}
M_h u(x)=\sum_{j} u(h x_j) \left(
   \frac{x_{j+1}-x/h}{x_{j+1}-x_{j}}\z_j\Big(\frac{x}{h} -x_j\Big)
+ \frac{x/h-x_{j-1}}{x_{j}-x_{j-1}} \,\z_{j-1} \Big(
\frac{x}{h} -x_{j-1} \Big)\right)
\end{split}
\end{equation*}
satisfies
$$
|M_h u(x)-u(x)|\leq c\, h^2\,  \Vert u''\Vert_{L_\infty(\R)} + \e \, |u(x)|\,,
$$
where the constant $c$ depends on the functions $\z_j$.

The outline of the paper is as follows.
In Section~\ref{sec_quasiuni} we consider an extension of 
the approximate quasi-interpolation to
scattered nodes close to a uniform grid. 
We construct the quasi-interpolant ${\mathbb M} u$ with gridded
centers and coefficients depending on scattered data and
obtain approximation estimates.
Further the results of 
some numerical experiments are presented which
confirm the predicted approximation orders.
In Section~\ref{S:DUE} we show that an approximate partition of unity
can be obtained from a given system of rapidly decaying approximating functions if
these functions are multiplied 
by polynomials. 
Using the approximate partition of unity,
one can construct 
approximate quasi-interpolants of high order
approximation rate up to some prescribed saturation error.
This is the topic of Section~\ref{S:DUE1}.
Section~\ref{S:CIN} contains 
an application 
to the cubature of convolution integral operators.
A construction of the approximate partition of unity
for the case of Gaussians and some numerical 
examples are given in Section~\ref{S:QUA}.

\paragrafo{Quasi-interpolants with gridded centers} \label{sec_quasiuni}

Here we give a simple extension of the quasi-interpolation operator on uniform grids
\eqref{1.1} to a quasi-interpolant,
which uses the values $u(\bx_\bj)$ on 
a set of scattered nodes $\bX = \{\bx_\bj\} \subset \R^n$ 
if it is close to a uniform grid.
Precisely we suppose 

\begin{cond} \label{cond1}
{\em There exist $h>0$ and $\k_{1} >0$ 
such that
for any $\bj\in \Z^{n}$ the ball $B(h \bj,h \k_{1})$ centered at
$h \bj$ with radius $h \k_{1}$ contains
nodes of $\bX$.
}
\end{cond}

\begin{defi} \label{defi1}
{\em
Let $\bx_{\bj}
\in \bX$. The collection of $\ds m_N=\frac{(N-1+n)!}{n!(N-1)!}-1$ nodes $\bx_{\bk}
\in \bX$
will be called the
{\it star} of $\bx_\bj$ and denoted by $\st(\bx_\bj)$ if
the Vandermonde matrix 
\begin{equation}\label{defVander}
V_{\bj,h}=
\Big\lbrace \Big( \frac{\bx_{\bk}- \bx_\bj}{h}\Big)^\ba
\Big\rbrace, \; |\ba|=1,...,N-1 ,
\end{equation}
 is not singular.
The union of the node $\bx_\bj$ and its star
is denoted by
$\ST(\bx_\bj) = \bx_\bj \cup \st(\bx_\bj)$.
}
\end{defi}

\begin{cond} \label{cond2}
{\em Assume Condition \ref{cond1} and denote
by $\widetilde \bx_{\bj} \in \bX$  
the node closest to $h\, \bj$.
There exists $ \k_{2}>0$ 
such that for any $\bj\in \Z^{n}$
\begin{itemize}
\item[(a)]
$\st(\widetilde \bx_\bj) \subset B(\widetilde \bx_{\bj},h\, \k_{2})$ 
with
$|\det V_{\bj,h} |\geq c > 0$;
\item[(b)] $ \bigcup \limits_{\bj \in \Z^n}  \ST(\widetilde \bx_\bj)= \bX$.
\end{itemize}
}
\end{cond}

\subparagrafo{Error estimate}
To formulate our first result we denote by
$\{b_{\ba,\bk}^{(\bj)}\}$, $\bx_{\bk} \in \st(\widetilde \bx_\bj)$, $|\ba|=1,...,N-1$,
the elements of the
inverse matrix of
$V_{\bj,h}$ and define the functional
\[
\Lambda_\bj(u) = u(\widetilde \bx_\bj)
\Big(1-\sum_{|\ba|=1}^{N-1}\Big(\bj-\frac{\widetilde \bx_\bj}{h}\Big)^\ba
\sum_{\bx_{\bk} \in \st(\widetilde \bx_\bj)}  b_{\ba,\bk}^{(\bj)} \,\Big)
+ \sum_{\bx_{\bk} \in \st(\widetilde \bx_\bj)} u(\bx_\bk) 
\sum_{|\ba|=1}^{N-1}
b_{\ba,\bk}^{(\bj)}\Big( \bj-\frac{\widetilde \bx_\bj}{h}\Big)^\ba \, .
\]

\begin{thm}\label{teorema0}
Suppose that for some  $K>n$ and  the smallest integer   $n_0 >n/2$
the function $\h(\bx)$, $\bx \in \R^n$,  satisfies the conditions
$(1+|\bx|)^K |\de^\bb \h(\bx)| \le C_{\bb}$  
for all $0 \le |\bb| \le n_0$, 
and $\de^\ba(\Ff \h -1)(\boldsymbol{0}) =0$,
$0 \le |\ba| < N$. 
If
the set of nodes $\bX$ satisfies Conditions {\em \ref{cond1}} and 
{\em \ref{cond2}},
then for any $\e>0$ 
there exists $D$ such that
the quasi-interpolant
\begin{equation}\label{quasiuniform}
{\mathbb M}\, u(\bx)= D^{-n/2} \sum_{\bj \in \Z^n} \Lambda_\bj(u)  
\h\Big( \frac{\bx-h\, \bj}{h \dd}\Big)
\end{equation}
approximates any $u \in  W_\infty^N(\R^n)$  with
\begin{equation}\label{est2bis}
|{\mathbb M} u(\bx)-u(\bx)|\leq c_{N,\h,D}\,\, h^{N} 
\sup_{\R^n}|\nabla_{N}u|+
\e \sum_{k=0}^{N-1}|\nabla_k u(\bx)| (\dd h)^{k} \, ,
\end{equation}
where $c_{N,\h,D}$ does not depend on $u$ and $h$.
\end{thm}
\textit{Proof.}
For given $u \in  W_\infty^N(\R^n)$ we consider
the quasi-interpolant \eqref{1.1} on the uniform grid $\{h \bj \}$
\[
\M\,u(\bx)=D^{-n/2} \sum_{\bj \in \Z^n} u(h\bj)\,
\h{\Big(\frac{\bx -h \bj}{h \sqrt{D}}\Big)} \, ,
\]
with $h$ given by Condition \ref{cond1}.
It was proved in \cite{MS4}
that under 
the  decay and moment conditions on $\h$, formulated in the statement
of the theorem, 
$\M\,u$ can be represented as
\[
\M\,u(\bx) = u(\bx) +  
\sum _{ |\ba| = \boldsymbol{0}} ^{N-1}
\Big(\frac{\dd h }{2 \pi i } \Big)^{|\ba|} 
      \frac{\de ^{\ba} u(\bx)}{\ba!} 
\sum_{\bk \in {\Z}^n\setminus \{\boldsymbol{0}\}} 
\de ^\ba \Ff \h (\dd \bk) \,
\ee^{2 \pi i (\bk,\bx)}
+ U_N ({\mathbf x})
\]
with a function $U_N$  bounded by 
\[
|U_N ({\mathbf x})| \le C_1\,\, (\dd h)^{N}\sup_{\R^n}|\nabla_{N}u|
\]
and a constant $C_1$ depending only on $\h$.
Moreover, the sequences
$\{ \de ^{\ba} \Ff 
\h(\dd \; \cdot)\} \in l_1 ({\Z}^n)$
and
\[
\sum_{\bk \in {\Z}^n \setminus \{\boldsymbol{0}\}}
|\de ^{\ba} \Ff  \h( \dd \bk)| \to 0 \quad \mbox{as} \quad D \to
\infty \, . 
\] 
Hence, we can find $D$ such that $\M u$ satisfies the inequality
\[
|\M u(\bx)-u(\bx)|\leq |U_N| + 
\e \sum_{k=0}^{N-1}|\nabla_k u(\bx)| (\dd h)^{k} \, .
\]
It remains to estimate
$|\M u - {\mathbb M} u|$.
Recall the 
Taylor expansion of $u$ around $\by \in \R^n$
\begin{equation}\label{Taylor}
u(\bx)=\sum_{|\ba| =0}^{N-1} \frac{\de^\ba u(\by)}{\ba!}
(\bx-\by)^\ba+
R_N(\by,\bx) 
\end{equation}
with the remainder  satisfying
\begin{equation}\label{R1j}
|R_N(\by,\bx)|\leq c_N |\bx-\by|^N \,
\sup_{B(\by,|\bx-\by|)}\vert \nabla_N u\vert \, .
\end{equation}
For $\bj \in \Z^n$ we choose $\widetilde \bx_\bj \in \bX$ and use
\eqref{Taylor} with $\by=\widetilde \bx_\bj$.
We split
\[
\M u(\bx) = M^{(1)}u(\bx)+D^{-n/2} \sum_{\bj \in \Z^n} 
R_N(\widetilde \bx_\bj,h \bj ) \, \h\Big( \frac{\bx-h\, \bj}{h \dd}\Big) 
\]
with
\begin{equation}\label{quam2}
M^{(1)}u(\bx)= D^{-n/2} \sum_{\bj \in \Z^n} \sum_{|\ba| =0}^{N-1}
\frac{\de^{\ba} u(\widetilde \bx_{\bj})}{\ba!}\, 
(h \bj-\widetilde \bx_{\bj})^{\ba}
\, \h\Big( \frac{\bx-h\, \bj}{h \dd}\Big) \, .
\end{equation}
Because of $|h \bj-\widetilde \bx_{\bj}| \le \kappa_1 h$ for any $\bj$
we derive from \eqref{R1j}
\begin{equation}\label{est2}
|M^{(1)}u(\bx)-\M u(\bx)|\leq c_N \, (\kappa_1 h)^{N} 
D^{-n/2} \sum_{\bj \in \Z^n} \Big|\h\Big( \frac{\bx-h\, \bj}{h \dd}\Big)\Big|
\, \sup_{B(\bx,h \k_{1})} |\nabla_{N}u| \, .
\end{equation}
The next step is to approximate $\de^\ba u(\widetilde \bx_\bj),\> 1\leq|\ba|<
N$, by a linear combination of $u(\bx_\bk)$, 
$\bx_{\bk} \in \st(\widetilde \bx_\bj)$.
Let $\{a_{\ba}^{(\bj)}\}_{1\leq |\ba|\le N-1}$ 
be the unique solution of the linear
system with $m_N$ unknowns
\begin{equation}\label{sistem1}
\sum_{|\ba| =1}^{N-1} \frac{a_{\ba}^{(\bj)}}{\ba!}\, 
(\bx_\bk-\widetilde\bx_\bj)^\ba=
u(\bx_\bk)-u(\widetilde\bx_\bj)\, , 
\quad \bx_{\bk} \in \st(\widetilde\bx_\bj) \, .
\end{equation}
From (\ref{Taylor}) and (\ref{sistem1}) follows that
\[
\sum_{|\ba| =1}^{N-1}\frac{h^{|\ba|}}{\ba!}
(a_{\ba}^{(\bj)}-\de^\ba u(\widetilde \bx_\bj))\,
\Big(\frac{\bx_{\bk}-\widetilde \bx_\bj}{h}\Big)^{\ba}=
R_N(\widetilde \bx_\bj,\bx_{\bk}) \, .
\]
By Condition \ref{cond2}(a)
the norms of $V_{\bj,h}^{-1}$
are bounded uniformly in $\bj$, this leads together with (\ref{R1j})
to the inequality
\begin{equation}\label{smdeu}
\frac{|a_{\ba}^{(\bj)}-\de^{\ba}u(\widetilde \bx_\bj)|}{ \ba!} 
\leq  C_2\, h^{N-|\ba|}
\sup_{B(\widetilde \bx_\bj,h \k_2 )}\vert\nabla_N u\vert.
\end{equation}
Hence, if we replace  the derivatives 
$\de^\ba u(\widetilde\bx_\bj)$ in \eqref{quam2} by
$a_{\ba}^{(\bj)}$, then we get the sum
\begin{equation*}
D^{-n/2} \sum_{\bj \in \Z^n}
\Big(u(\widetilde \bx_\bj) + \sum_{|\ba|=1}^{N-1}
\frac{a_{\ba}^{(\bj)} }{\ba!}\, (h \bj-\widetilde \bx_{\bj})^{\ba} \Big) \,
\h\Big( \frac{\bx-h\, \bj}{h \dd}\Big),
\end{equation*}
which in view of
\[
a_{\ba}^{(\bj)}=\frac{\ba!}{h^{|\ba|}} 
\sum_{\bx_\bk\in \st(\widetilde \bx_\bj)}
b_{\ba,\bk}^{(\bj)} \, (u(\bx_\bk)-u(\widetilde \bx_\bj))
\]
coincides with \eqref{quasiuniform}.
Moreover,
\begin{equation}\label{est3}
|{\mathbb M}\, u(\bx)- M^{(1)}u(\bx)|\leq C_2 h^{N} 
\sum_{|\ba|=1}^{N-1} \k_1^{|\ba|} \, 
D^{-n/2}  \sum_{\bj \in \Z^n} \Big|\h\Big( \frac{\bx-h\, \bj}{h \dd}\Big)\Big|
  \sup_{B(\bx,h \k_{2})}\vert\nabla_N u\vert \, .
\end{equation}
Now the inequality
\[
\sup_{\R^n} D^{-n/2}  
\sum_{\bj \in \Z^n} \Big|\h\Big( \frac{\bx- \bj}{\dd}\Big)\Big|
\le C_3
\]
for all $D \ge D_0 > 0$ implies that \eqref{est2} and \eqref{est3}
lead to
\[
|\M u(\bx) - {\mathbb M}\,u (\bx)|
\le C_4 h^{N}   \sup_{B(\bx,h \k_{2})}\vert\nabla_N u\vert \, ,
\]
which proves \eqref{est2bis}.
\hfill $\blacksquare$

\subparagrafo{Numerical Experiments with Quasi-interpolants}\label{33}
The behavior of the quasi-interpolant ${\mathbb M}u$ 
was tested by one- and two-dimensional experiments. 
In all cases the
scattered grid is chosen such that
any ball $B(h \bj, h/2)$, $\bj \in \Z^n$,  contains one randomly
chosen node $\bx_\bj$, thus $\widetilde \bx_\bj= \bx_\bj$.
All the computations were carried out with MATHEMATICA\textregistered.

\begin{figure}[h] \begin{center}
 \epsfig{file=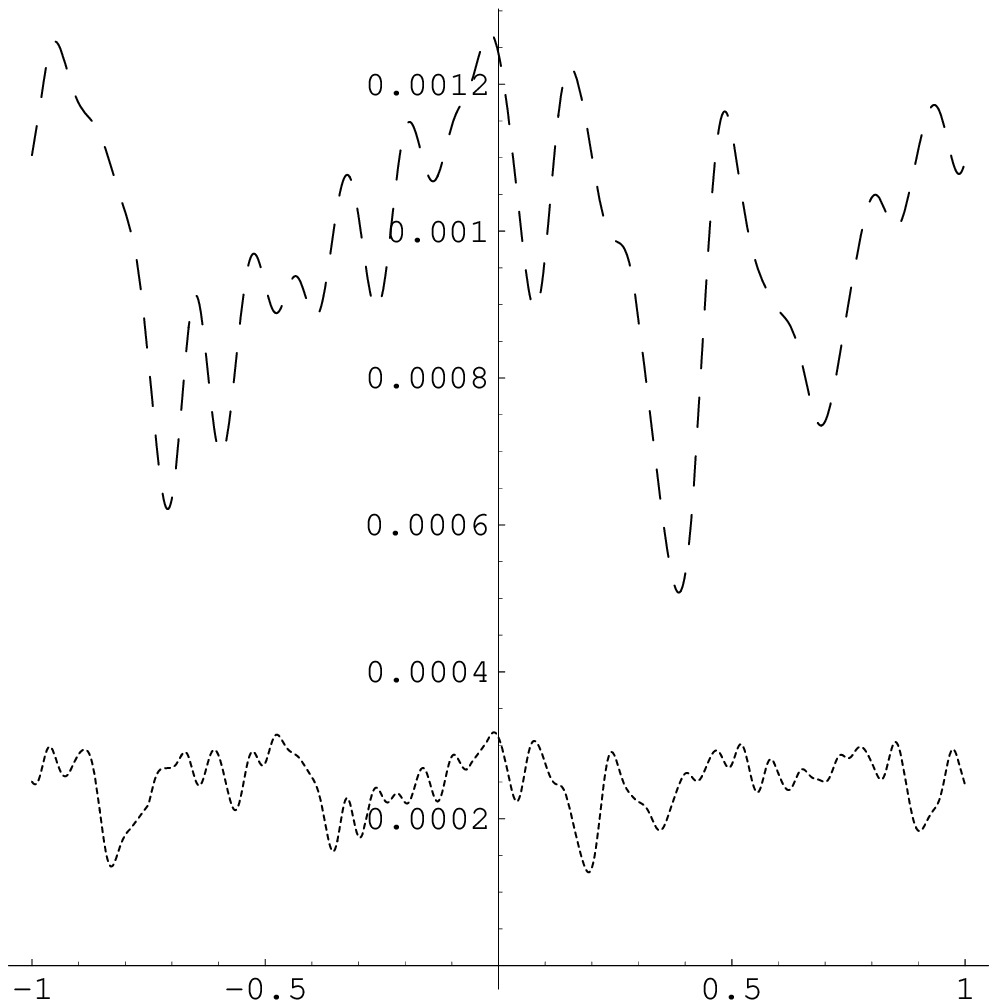,width=3in,height=1.8in}
 \epsfig{file=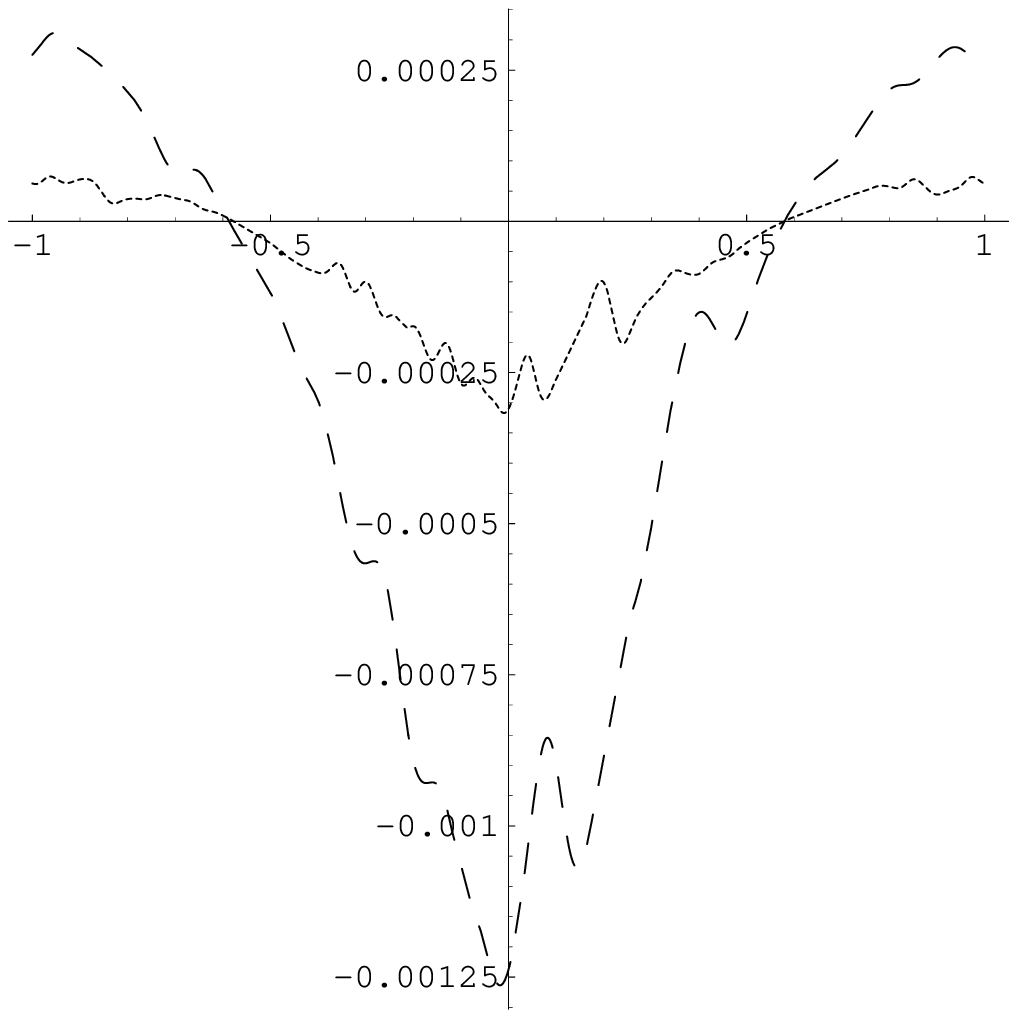,width=3in,height=1.8in}
 \caption{{\small The graphs of ${\mathbb M}\,u(x)-u(x)$  with
$D=2$, $N=2$, $\st(x_j)=\{x_{j+1}\}$, when $u(x)=x^2$ (on the left) and 
 $u(x)=1/(1+x^{2})$. Dashed and solid lines correspond to 
 $h=1/32$ and   $h=1/64$.}}
 \label{fig1}
 \end{center}
 \end{figure}

\begin{figure}[h] \begin{center}
 \epsfig{file=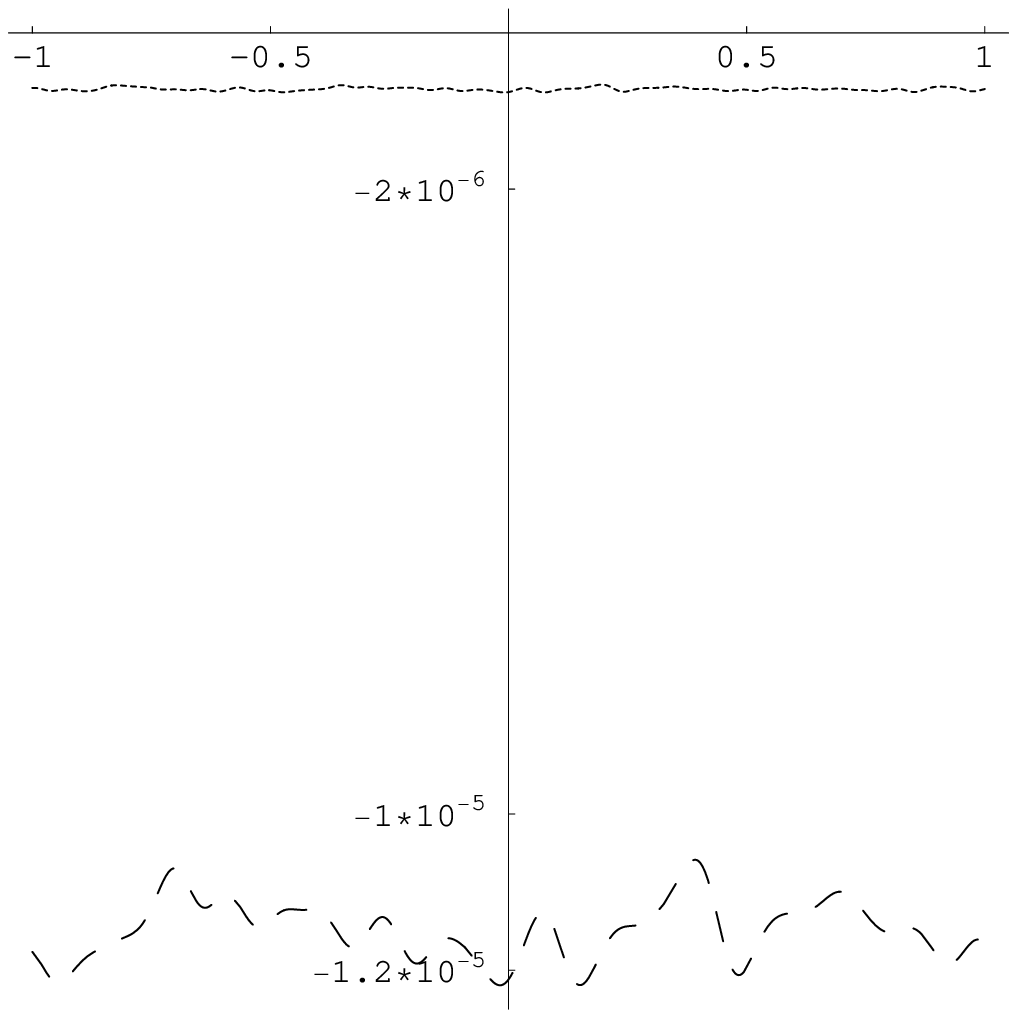,width=3in,height=1.8in}
 \epsfig{file=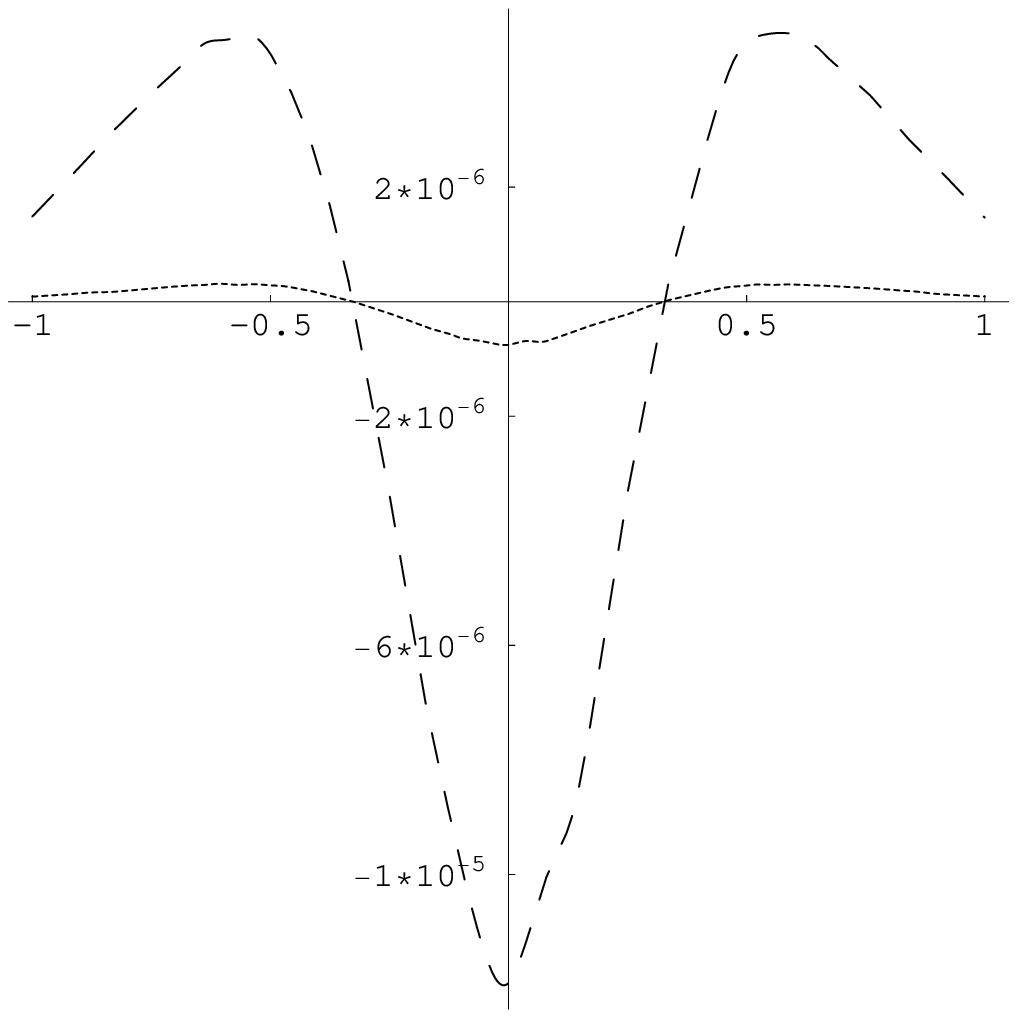,width=3in,height=1.8in}
 \caption{{\small The graphs of ${\mathbb M}\,u(x)-u(x)$  with
$D=4$, $N=4$, $\st(x_j)=\{x_{j-2}, x_{j-1}, x_{j+1}\}$,
 when $u(x)=x^4$ (on the left) and 
 $u(x)=1/(1+x^{2})$. Dashed and solid lines correspond to 
 $h=1/32$ and   $h=1/64$.}}
 \label{fig2}
 \end{center}
 \end{figure}

The following  figures show the graph of
${\mathbb M}u-u$ for different smooth functions $u$
using basis functions for second (Fig.~\ref{fig1}) 
and fourth (Fig.~\ref{fig2}) order of approximation 
with   $h=1/32$ (dashed line) and $h=1/64$ (solid line).
 

\begin{figure}[h] \begin{center}
 \epsfig{file=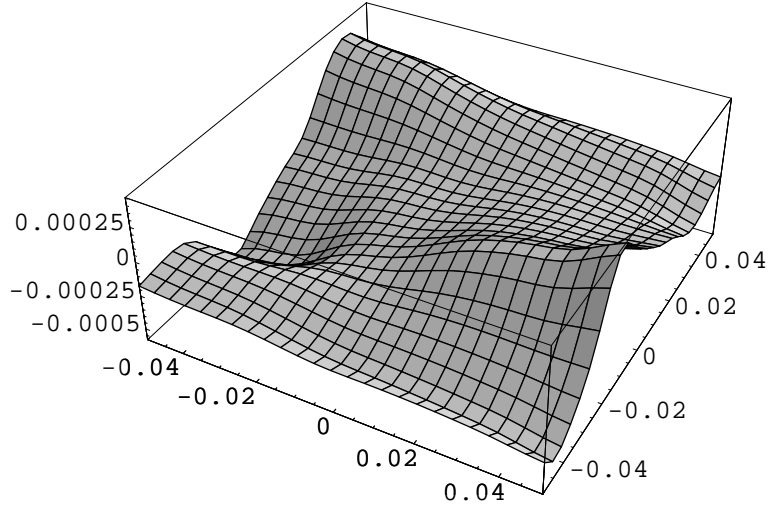}
 \caption{{\small The graph of ${\mathbb M} u-u$  with
$D=2$, $N=2$, $h=1/128$, when $u(\bx)=1/(1+|\bx|^{2})$, $\bx \in \R^2$.
}}
 \label{fig3}
 \end{center}
 \end{figure}

In Fig.~\ref{fig3} the difference ${\mathbb M}u-u$   
is plotted for the function $u(\bx)=(1+|\bx|^{2})^{-1}$, $\bx \in
\R^2$. 
In formula (\ref{quasiuniform}) we use 
$\h(\bx)=\pi^{-1}{\ee}^{-|\bx|^{2}}$, $D=2$, for which $N=2$, and
therefore the star $\st(\bx_\bj)$ contains $2$ nodes.
We have chosen $\st(\bx_{j_1,j_2})= \{\bx_{j_1+1,j_2},\bx_{j_1,j_2+1}\}$.
In Table  \ref{Table1} we give some computed values of 
${\mathbb M}\,u(\boldsymbol{0})-u(\boldsymbol{0})$ for $D=2$ and $D=4$
with different $h$, which 
confirm $h^2$-convergence of the two-dimensional quasi-interpolant.
\begin{table}[h] 
    \centering
    \begin{tabular}{|c|c|c|}\hline
	$h$ & $D=2$ & $D=4$ \\  \hline
	$2^{-4}$ & $-6.2 \cdot 10^{-3}$ & $-1.3 \cdot 10^{-2}$\\
	$2^{-5}$ & $-1.6 \cdot 
10^{-3}$ & $-3.3 \cdot 10^{-3}$\\
	$2^{-6}$ & $-3.9 \cdot 10^{-4}$ & $-8.3 \cdot 10^{-4}$\\
	$2^{-7}$ & $-9.8 \cdot 
10^{-5}$ & $-2.1 \cdot 10^{-4}$\\
	$2^{-8}$ & $-2.4 \cdot 10^{-5}$ & $-5.2 \cdot 10^{-5}$\\
	\hline
\end{tabular}	
\caption{Values of ${\mathbb M}\,u(\boldsymbol{0})-u(\boldsymbol{0})$}   
\label{Table1}
\end{table}

\paragrafo{Approximate partition of unity}\label{S:DUE}
\setcounter{equation}{0}
In the following two sections we consider
irregularly distributed nodes.
First we show that an approximate partition of unity
can be obtained from a given system of approximating functions
centered at the scattered nodes
if
these functions are multiplied 
by polynomials. 
We are mainly interested in rapidly decaying basis functions which are
supported on the whole space.
But we start with the simpler case of compactly supported
basis functions.  

\subparagrafo{Basis functions with compact support}\label{SUB:COM}

\begin{lemma}\label{lemma1}
Let $\{B(\bx_\bj,h_\bj)\}_{\bj \geq 0}$ be an open locally finite covering
of $\R^n$ by balls centered in $\bx_\bj$ and radii $h_\bj$. Suppose
that the multiplicity
of this covering does not exceed a positive constant $\m_n$ and that
there are positive
constants $c_1$ and $c_2$ satisfying
\begin{equation}\label{hjm}
c_1 h_\bm \leq h_\bj \leq c_2 h_\bm
\end{equation}
provided the balls $B(\bx_\bj,h_\bj)$ and $B(\bx_\bm,h_\bm)$ have
common points. Furthermore, let $\{\hj\}$ be a
bounded sequence of continuous functions on $\R^n$ such that $\supp
\h_\bj \subset B(\bx_\bj,h_\bj)$.
We assume that the functions
$\R^n \ni \by \to \h_\bj(h_\bj\, \by)$ are continuous uniformly with
respect to $\bj$
and
\begin{equation}\label{lowbound}
s(\bx):=\sum_{\bj} \hj(\bx)\geq c \>\>\>{\rm on}\>\>\R^n
\end{equation}
where $c$ is a positive constant.
Then for any $\e>0$ there exists a sequence of polynomials $\{\Pol_\bj\}$
with the following properties:
\begin{itemize}
\item[\rm (i)] the degrees of all $\Pol_\bj$ are bounded (they depend on the
least majorant
of the continuity modulae of $\h_\bj$ and the constants $\e$, $c$,
$c_1$, $c_2$, $\m_n$);
\item[\rm (ii)] there is such a constant $c_0$ that
$|\Pol_\bj|<c_0 $ on $B(\bx_\bj,h_\bj)$;
\item[\rm (iii)] the function
\begin{equation}\label{th3}
\Th:=\sum_{\bj} \Pol_\bj \,\h_\bj \end{equation}
   satisfies
   \begin{equation}\label{tm1}
|\Th(\bx)-1|<\e \quad \mbox{ for all } \,  \bx \in \R^n \,.
\end{equation}
\end{itemize}
\end{lemma}
\textit{Proof.} 
Since the functions $B(\bx_\bj,1) \ni \by \to s(h_\bj\, \by)$
are continuous uniformly with
respect to $\bj$, 
for an arbitrary positive $\delta$ 
there exist polynomials $\Pol_{\bj}$
subject to
\begin{equation*}
\Big|\Pol_{\bj}(\bx)- \frac{1}{s(\bx)}\Big|<\delta \>\>{\rm on}\>\> 
B(\bx_\bj,h_\bj) \, ,
\end{equation*}
and the degree of $\Pol_{\bj}$, $\deg \Pol_{\bj}$,
is independent of $\bj$.
Letting $\delta =\e \, (\m_n \|\h\|_{L_\infty})^{-1}$
we obtain
\begin{equation}\label{infi_ex1}
\Big|\h_\bj(\bx) \Big(\Pol_{\bj}(\bx) - \frac{1}{s(\bx)}\Big) \Big|
\le \frac{\e}{\m_n} \, . 
\end{equation}
Then
\begin{equation}\label{tm1ter}
\sup_{\R^n} \sum_{\bj}\Big|\h_\bj(\bx) \Big(\Pol_{\bj}(\bx) -
\frac{1}{s(\bx)}\Big) \Big|
\le \e \, ,
\end{equation}
since at most $\m_n$ terms of this sum are different from zero.
But
\begin{align*}
\sum_{\bj}& \h_\bj(\bx) \Big(\Pol_{\bj}(\bx) -\frac{1}{s(\bx)}
\Big) 
= \sum_{\bj} \h_\bj(\bx) \Pol_{\bj} (\bx) - \frac{1}{s(\bx)}
\sum_{\bj} \h_\bj(\bx)   
= \sum_{\bj} \h_\bj(\bx) \Pol_{\bj} (\bx) - 1 \, ,
\end{align*}
which proves (\ref{tm1}). \hfill $\blacksquare$
\bigskip

\begin{rem} {\em 
Let the functions $\{\hj\}$ in Lemma
\ref{lemma1} satisfy
the additional hypothesis $\h_\bj \in C^k(\R^n)$. Then one can find a
sequence of polynomials $\{\Pol_{\bj}\}$ of degrees
$L_{\bj}$ such that
$$
\sup_{B(\bx_{\bj},h_{\bj})}\Big|\Pol_{\bj}(\bx)- \frac{1}{s(\bx)}\Big|
\leq
C(k)\, \frac{h^k_\bj}{L_{\bj}^k} \sup_{B(\bx_\bj,h_\bj)}\vert
\nabla_k {s(\bx)}\vert
$$
(see, \textit{e.g.}, \cite{Mha}). This
shows that it suffices to take
polynomials $\Pol_{\bj}$ with $\deg \Pol_{\bj} > c(k)\, \e^{-1/k}$
in order to achieve the error $\e$ in \eqref{tm1ter}.
}
\end{rem}

\subparagrafo{Basis functions with noncompact support}\label{SUB:NOCOM}
Here we consider  approximating functions  supported
on the whole $\R^n$.
We suppose that the functions $\h_\bj$ are scaled translates 
\[
\h_\bj(\bx) = \h \Big(\frac{\bx -\bx_\bj}{h_\bj}\Big) 
\]
of a sufficiently smooth function $\h$ with rapid decay.
\begin{lemma}\label{lemma211}
For any $\e > 0$ there exists $L_{\e}$ and 
polynomials $\{\Pol_\bj\}$ of degree $\deg \Pol_\bj \le L_{\e}$
such that the function $\Th$ defined by {\em(\ref{th3})}
  satisfies
{\em(\ref{tm1})}
under the following assumptions on $\eta$, the nodes
$\{\bx_\bj\}$ and the scaling parameters $\{h_\bj\}$:
\begin{itemize}
\item[\rm 1.] There exists $K>0$ such that
\begin{equation}\label{distri_ex11}
c_K:=\Big\|\sum_{\bj} \big(1+ h_\bj^{-1} |\,{\cdot-\bx_\bj}|\big)^{-K} 
\Big\|_{L_\infty}  < \infty \, .
\end{equation}
\item[\rm 2.] There exists $p>0$ such that
\begin{equation}\label{decay11}
\Big\|\big(1+  |\cdot | \big)^{K} \ee^{p^2|\cdot|^2}\, 
\h \Big\|_{L_\infty}\, , \;
\Big\|\big(1+  |\cdot | \big)^{K} \,\nabla  \h 
\Big\|_{L_\infty}  \le c_p < \infty \, .
\end{equation}
\item[\rm 3.] There exists $C\ge 1$ such that all indices $\bj$, $\bm$
\begin{equation}\label{h_bound11}
\frac{h_\bj} {h_\bm}  \le C \, .
\end{equation}
\item[\rm 4.] \eqref{lowbound} is valid.
\end{itemize}
\end{lemma}
\textit{Proof.} 
From 
\eqref{distri_ex11} and \eqref{decay11}
the sum 
$
s(\bx)=\sum_{\bj} \h_\bj(\bx) 
$
converges absolutely for any $\bx$ to a positive, smooth and bounded function $s$.
Suppose that we have shown that
for  any $\e>0$ and all indices $\bj$ there exist polynomials
$\Pol_{\bj}$ such that
\begin{equation}\label{infi_ex11}
\Big|\h_\bj(\bx) \Big(\Pol_{\bj}(\bx) - \frac{1}{s(\bx)}\Big) \Big|
\le \frac{\e}{c_K} \, \big(1+ {h_\bj}^{-1} |{\bx-\bx_\bj}| \big)^{-K} \, ,
\end{equation}
($c_K$ is defined in  \eqref{distri_ex11})
and $\deg \Pol_{\bj} \le L_{\e}$.
Then
\begin{equation*}
\sup_{\R^n} \sum_{\bj}\Big|\h_\bj(\bx) \Big(\Pol_{\bj}(\bx) -
\frac{1}{s(\bx)}\Big) \Big|
\le \e  \, ,
\end{equation*}
and as in the proof of Lemma \ref{lemma1} we conclude
\[
\sup_{\R^n}  \Big|\sum_{\bj} \h_\bj(\bx) \Pol_{\bj} (\bx) - 1\Big| 
\le  \e  \, .
\]  
To establish (\ref{infi_ex11}) we use
a result on weighted polynomial approximation 
from
\cite{DzT}, which will be stated in a simplified form.
Introduce the weight  $w_p(\bx)=\ee^{-p^2|\bx|^2}$, $p > 0$, and 
consider the 
best weighted polynomial approximation of
a function $g$ given on $\R^n$
\begin{equation*}
E_N(g )_{w_p,\infty}  := 
\inf_{ \Pol \in \Pi_N} \sup_{\R^n} 
|w_p(\bx)(g(\bx) - \Pol(\bx))| \, ,
\end{equation*}
where $\Pi_N$ denotes the set of polynomials, which are of degree at
most $N$ in each variable $x_1,\ldots,x_n$.             
Then for  $g \in W_\infty^1(\R^n)$ 
\begin{equation} \label{westim11}
E_N(g )_{w_p,\infty} \le C  \,
\frac{\|\nabla g \|_{L_\infty}}{N^{1/2}} \, .
\end{equation} 
Let us fix an index $\bj$
and make the change of variables $\by={h_\bj}^{-1} ({\bx-\bx_\bj})$.
Then  (\ref{infi_ex11}) is proved if we show that there
exists a polynomial $\Pol_\bj$ such that for all $\by \in \R^n$ 
\begin{equation} \label{polest11}
\Big| \h (\by)  \Big(\Pol_{\bj}(\by) - 
\frac{1}{\tilde s(\by)}\Big) \Big|
\le \frac{\e}{c_K} \big(1+  |\by| \big)^{-K} 
\end{equation}
with $\tilde s(\by) = s(h_\bj \by +\bx_\bj)$.
Since
$\tilde s ^{-1} \in W^1_\infty(\R^{n})$
according to (\ref{westim11}) we can find a 
polynomial $\Pol_{\bj}$ satisfying
\[
\sup_{\R^n}\Big|\Pol_{\bj}(\by)- \frac{1}{\tilde s(\by)}\Big|
\ee^{-p^2|\by|^2} < \frac{\e}{c_p \, c_K} 
\]
with the constant $c_p$ in the decay condition
(\ref{decay11}). Now (\ref{polest11}) follows immediately from
\[
   |\h (\by)| \big(1+  |\by| \big)^{K} \le c_p \, \ee^{-p^2|\by|^2} \, .
\]
From (\ref{westim11}) we see 
that $\deg \Pol_{\bj}$ depends on the norm of the gradient
\begin{equation*} \label{gradient11}
\begin{split}
\sup_{\R^n}
\Big|  \nabla \frac{1}{\tilde s(\by)} \Big|
&= h_\bj  \sup_{\R^n}
\Big| \nabla \frac{1}{s(\bx)} \Big| 
\le  \sup_{\R^n}
\frac{1}{(s(\bx))^{2}} 
\sum_{\bm} \frac{h_\bj} {h_\bm}
\Big| \nabla \h \Big(\frac{\bx -\bx_\bm}{h_\bm}\Big) \Big| \, ,
\end{split}
 \end{equation*}
which is bounded uniformly in $\bj$ 
in view of 
\eqref{lowbound}, \eqref{decay11}, and \eqref{h_bound11}.
\hfill $\blacksquare$

\paragrafo{Quasi-interpolants of a general form}\label{S:DUE1}
\setcounter{equation}{0}
In this section we study the approximation of functions 
$u\in W_\infty^N(\R^n)$ by the
quasi-interpolant (\ref{e:nonuniform}). 
We will show that within the class of generating functions
of the form polynomial times compactly supported or rapidly decaying
generating function it suffices to have an approximate partition of unity
in order to construct approximate quasi-interpolants of high order
accuracy up to some prescribed saturation error.

Let us assume the following hypothesis concerning the grid $\{\bx_\bj\}$:

\begin{cond} \label{cond3}
{\em
For any $\bx_\bj$ there exists a ball
$B(\bx_\bj,h_\bj)$
which contains
$m_N$
nodes $\bx_{\bk} \in \st(\bx_\bj)$ with
\begin{equation}\label{Vandermonde1}
|\det V_{\bj, h_\bj} | = \Big| \det \Big\lbrace 
\Big( \frac{\bx_{\bk}-\bx_\bj}{h_\bj}\Big)^\ba \Big\rbrace_{|\ba|=1}^{N-1}
 \Big| \geq c  \, ,
\end{equation}
(see Definition \ref{defVander}), with $c>0$ not depending on $\bx_\bj$.
}
\end{cond}

\subparagrafo{Compactly supported basis functions} \label{31}
\begin{thm}\label{teorema1}
Suppose that the function system $\{\h_\bj\}$ satisfies 
the conditions of Lemma {\em \ref{lemma1}},   
let $u\in W_\infty^N(\R^n)$ and $\e > 0$ arbitrary. 
There exist polynomials $\Pol_{\bj,\bk}$, independent
on $u$,
whose degrees are uniformly bounded, such that the quasi-interpolant
\begin{equation}\label{mainM}
Mu(\bx)=\sum_{\bk} u(\bx_\bk) \sum_{\ST(\bx_\bj)\ni \bx_\bk}
\Pol_{\bj,\bk}(\bx)\h_\bj(\bx)
\end{equation}
satisfies the estimate
\begin{equation}\label{mainR}
|Mu(\bx)-u(\bx)|\leq C h_\bm^N\,\sup_{{B}(\bx_\bm,\l\,h_\bm)} |\nabla_N
\, u|\,  +\e\,|u(\bx)|,
\end{equation}
where $\bx_\bm$ is an arbitrary node and $\bx$ is any point
of the ball  $B(\bx_\bm,h_\bm)$.
By $\l$ we denote a constant greater than $1$ which depends on $c_1$
and $c_2$ in \eqref{hjm}.
The constant $C$ does not depend on $h_\bm$, $\bm$ and $\e$.
\end{thm}
\textit{Proof.} For given $\e$ we choose polynomials
$\Pol_\bj(\bx)$ such that the function (\ref{th3}) satisfies
\[
|\Th(\bx)-1|<\e \quad \mbox{ for all } \,  \bx \in \R^n   \, ,
\]
and introduce the auxiliary quasi-interpolant
\begin{equation}\label{Mtilde}
{ M^{(1)}}u(\bx)=\sum_{\bj} \Big(\sum_{|\ba|=0}^{N-1} \frac{\de^\ba
u(\bx_\bj)}{\ba!}(\bx-\bx_\bj)^\ba \Big) \Pol_\bj(\bx)\hj \asb \, .
\end{equation}
Using the Taylor expansion \eqref{Taylor} with $\by=\bx_\bj$
we write $M^{(1)}u(\bx)$
as
\[
{M}^{(1)} u(\bx)= u(\bx) {\Th}(\bx)
   -\sum_{\bj} R_N(\bx_\bj,\bx)\Pol_\bj(\bx) \hj \asb \,,
\]
which gives
\begin{equation*}
|M^{(1)}u(\bx)-u(\bx)|\leq  \sum_{\bj} |R_N(\bx_\bj,\bx)\Pol_\bj(\bx) \hj
\asb|+|u(\bx)|\,
|\Th(\bx)-1|\, .
\end{equation*}
This, together with the estimate for the remainder (\ref{R1j}), 
shows that for $\bx \in B(\bx_\bm,h_\bm)$
\begin{equation}\label{est1}
|M^{(1)}u(\bx)-u(\bx)|\leq C_1 h_\bm^N\,\sup_{B(\bx_\bm,\l\,h_\bm)}
|\nabla_N \, u|+
\e|u(\bx)| \, ,
\end{equation}
where the ball $B(\bx_\bm,\l\,h_\bm)$ contains all balls
$B(\bx_\bj,\,h_\bj)$
such that $B(\bx_\bj,\,h_\bj)$ and
$B(\bx_\bm,\,h_\bm)$ intersect.

Similar to the proof of Theorem \ref{teorema0}
we approximate in ${M}^{(1)} u$ 
the values of the derivatives  $\de^\ba u(\bx_\bj)$
by a linear combination of $u(\bx_\bk)$, 
where $\bx_{\bk} \in \st(\bx_\bj)$.
The solution of
\begin{equation*}\label{sistema}
\sum_{|\ba|=1}^{N-1} \frac{a_{\ba}^{(\bj)}}{\ba!}\, (\bx_\bk-\bx_\bj)^\ba=
u(\bx_\bk)-u(\bx_\bj)\, , \quad \bx_{\bk} \in \st(\bx_\bj) \, , 
\end{equation*}
is given by
\[
a_{\ba}^{(\bj)}=\frac{\ba!}{h_\bj^{|\ba|}} 
\sum_{\bx_\bk\in \st(\bx_\bj)}
b_{\ba,\bk}^{(\bj)} \, (u(\bx_\bk)-u(\bx_\bj)),
\]
where 
$\{b_{\ba,\bk}^{(\bj)}\}$ are the elements of the
inverse of $V_{\bj,h_\bj}$. 
Replacing the derivatives $\{\de^\ba u(\bx_\bj)\}$ in (\ref{Mtilde}) by
$\{a_{\ba}^{(\bj)}\}$ gives the quasi-interpolant
\begin{align*}
Mu(\bx)= & \, \sum_{\bj} \Big\{u(\bx_\bj)
\Big(1-\sum_{\bx_\bk\in \st(\bx_\bj)}\sum_{|\ba| =1}^{N-1}
b_{\ba,\bk}^{(\bj)}
\Big(\frac{\bx-\bx_\bj}{h_\bj}\Big)^\ba
\Big) \\& \,\qquad 
+ \sum_{\bx_\bk\in \st(\bx_\bj)}u(\bx_\bk)
\sum_{|\ba|  =1}^{N-1}
b_{\ba,\bk}^{(\bj)}
\Big(\frac{\bx-\bx_\bj}{h_\bj}\Big)^\ba
\,\Big\}\,\Pol_\bj(\bx)\h_\bj(\bx)\\
= & \, \sum_{\bj} \sum_{\bx_\bk\in \ST(\bx_\bj)}u(\bx_\bk)\, 
\Pol_{\bj,\bk}(\bx)\,\h_\bj(\bx)
\end{align*}
which can be rewritten as the quasi-interpolant (\ref{mainM}).
By (\ref{Taylor}) we obtain again
\[
\sum_{|\ba| =1}^{N-1}\frac{h_\bj^{|\ba|}}{\ba!}
(a_{\ba}^{(\bj)}-\de^\ba u(\bx_\bj))\,
\Big(\frac{\bx_{\bk}-\bx_\bj}{h_\bj}\Big)^{\ba}=
R_N(\bx_\bj,\bx_{\bk}) \, ,
\]
hence the boundedness of $\|V_{\bj,h_\bj}^{-1}\|$ from  Condition \ref{cond3}
and the estimate of the remainder (\ref{R1j})
imply  
\begin{equation*}
|a_{\ba}^{(\bj)}-\de^{\ba}u(\bx_\bj)|\leq \ba! \, C_2\, h_\bj^{N-|\ba|}
\sup_{B(\bx_\bj,h_\bj)}\vert\nabla_N u\vert.
\end{equation*}
Therefore we obtain the inequality
\[ 
|Mu(\bx)-{M}^{(1)}u(\bx)|\leq C_2 \sum_{\bj} h_\bj^{N}
\sup_{B(\bx_\bj,h_\bj)}\vert\nabla_N
u\vert\sum_{|\ba|=1}^{N-1} 
\Big|\frac{\bx-\bx_\bj}{h_\bj}\Big|^{|\ba|}\, 
|\Pol_\bj(\bx)\h_\bj(\bx)|
\]
and, for any $\bx\in B(\bx_\bm,h_\bm)$,
\[
|Mu(\bx)-{M}^{(1)}u(\bx)|\leq C_3 \, 
h_\bm^N\,\sup_{B(\bx_\bm,\l\,h_\bm)}\vert\nabla_N u\vert.
\]
This inequality  and (\ref{est1}) lead to (\ref{mainR}). \hfill $\blacksquare$

\subparagrafo{Quasi-interpolants with 
noncompactly supported basis functions}\label{32}
\begin{thm}\label{teorema1a}
Suppose that additionally to the
conditions of Lemma {\em \ref{lemma211}} 
the inequality
\begin{equation}\label{distri_exn}
\Big\|\sum_{\bj} \big(1+ h_\bj^{-1} |\,{\cdot-\bx_\bj}|\big)^{N-K} 
\Big\|_{L_\infty}  < \infty 
\end{equation}
is fulfilled, 
let $u\in W_\infty^N(\R^n)$ and $\e > 0$ arbitrary. 
There exist polynomials $\Pol_{\bj,\bk}$, independent
on $u$,
whose degrees are uniformly bounded, such that the quasi-interpolant
\begin{equation}\label{mainM1}
Mu(\bx)=\sum_{\bk} u(\bx_\bk) \sum_{\ST(\bx_\bj)\ni \bx_\bk}
\Pol_{\bj,\bk}\Big(\frac{\bx -\bx_\bj}{h_\bj}\Big) 
\h \Big(\frac{\bx -\bx_\bj}{h_\bj}\Big)
\end{equation}
satisfies the estimate
\begin{equation}\label{mainR1}
|Mu(\bx)-u(\bx)|\leq C \sup_{\bm}
h_\bm^N\,
\, \|\nabla_N \, u\|_{L_\infty}\,  +\e\,|u(\bx)| \, .
\end{equation}
The constant $C$ does not depend on $u$ and $\e$.
\end{thm}
\textit{Proof.} Analogously to \eqref{Mtilde} 
we introduce the quasi-interpolant 
\[
{ M^{(1)}}u(\bx)=\sum_{\bj} \Big(\sum_{|\ba| =0}^{N-1} \frac{\de^\ba
u(\bx_\bj)}{\ba!}(\bx-\bx_\bj)^\ba \Big) 
\Pol_\bj\Big(\frac{\bx -\bx_\bj}{h_\bj}\Big) 
\h \Big(\frac{\bx -\bx_\bj}{h_\bj}\Big)
\]
and obtain the estimate
\[
|M^{(1)}u(\bx)-u(\bx)| \leq \sum_{\bj} \Big|R_N(\bx_\bj,\bx)
\Pol_\bj \Big(\frac{\bx -\bx_\bj}{h_\bj}\Big) 
\h \Big(\frac{\bx -\bx_\bj}{h_\bj}\Big)\Big|
+|u(\bx)|\, |\Th(\bx)-1| \, .
\]
From \eqref{polest11} we have
\[
\Big|\Pol_\bj \Big(\frac{\bx -\bx_\bj}{h_\bj}\Big) 
\h \Big(\frac{\bx -\bx_\bj}{h_\bj}\Big)\Big|
\le \frac{1}{c} \Big|\h \Big(\frac{\bx -\bx_\bj}{h_\bj}\Big)\Big|
+ \frac{\e}{c_K}\Big(1+ \frac{|{\bx-\bx_\bj}|}{h_\bj}\Big)^{-K} 
\]
with the lower bound $c$ of $s(\bx)$ (see \eqref{lowbound}).
Together with \eqref{decay11} and \eqref{R1j} this provides
\begin{equation*}
\begin{split}
\Big|R_N(\bx_\bj,\bx)&
\Pol_\bj \Big(\frac{\bx -\bx_\bj}{h_\bj}\Big) 
\h \Big(\frac{\bx -\bx_\bj}{h_\bj}\Big)\Big| \\
&\le c_N h_\bj^N \|\nabla_N \, u\|_{L_\infty} 
\Big|\frac{\bx -\bx_\bj}{h_\bj}\Big|^{N}
\Big(\frac{c_{p}}{c} \, \ee^{-p^2 |\bx -\bx_\bj|^2/h_\bj^2}+ \frac{\e}{c_K}\Big)
\Big(1+\Big|\frac{\bx -\bx_\bj}{h_\bj}\Big|\Big)^{-K}
\end{split}
\end{equation*}
resulting in
\begin{equation*}
\begin{split}
|M^{(1)}&u(\bx)-u(\bx)| \leq  |u(\bx)|\, |\Th(\bx)-1| 
+ c_N \|\nabla_N \, u\|_{L_\infty}  \\
 & \times 
\Big(\frac{c_{p}}{c} \big\| \ee^{-p^2|\bx|^2} |\bx|^{N}\big\|_{L_\infty}
\sum_{\bj}h_\bj^N
\Big(1+\Big|\frac{\bx -\bx_\bj}{h_\bj}\Big|\Big)^{-K} 
\!\! +\frac{\e}{c_K}  
\sum_{\bj}h_\bj^N
\Big(1+\Big|\frac{\bx -\bx_\bj}{h_\bj}\Big|\Big)^{N-K}\Big) .
\end{split}
\end{equation*}
Now we can proceed 
as in the proof of Theorem \ref{teorema1}.
\hfill $\blacksquare$

\begin{rem} {\em 
Let for   fixed $\bx$ the parameter $\kappa_\bx$ be chosen such that
\[
\sum_{|\bx_\bj-\bx|> \kappa_\bx}   
\ee^{-p^2 |\bx -\bx_\bj|^2/h_\bj^2}
\Big|\frac{\bx -\bx_\bj}{h_\bj}\Big|^{N}
\Big(1+\Big|\frac{\bx -\bx_\bj}{h_\bj}\Big|\Big)^{-K} < \e \, .
\]
Then the estimate \eqref{mainR1} can be sharpened to
\begin{equation*}
|Mu(\bx)-u(\bx)|\leq C \max_{|\bx_\bj- \bx|\le \kappa_\bx}
h_\bj^N\,
\, \sup_{{B}(\bx,\kappa_\bx)} |\nabla_N
\, u|\,
+ \e\,(|u(\bx)| +\|\nabla_N \, u\|_{L_\infty}) \, .
\end{equation*}
}
\end{rem}

\paragrafo{Application to the computation of integral
operators}\label{S:CIN}
\setcounter{equation}{0}
Here we discuss a  direct application of the quasi-interpolation formula 
\eqref{mainM1} for the important example 
$ \h(\bx)= \ee^{-|\bx|^2}$.
Suppose that the density of the integral operator with radial
kernel
\begin{equation} \label{defint}
        \K u(\bx)=\int_{\R^{n}} g(|\bx - \by|) u(\by) \, d \by
\end{equation}
is approximated by
the quasi-interpolant 
\begin{equation}\label{mainM2}
M u(\bx)=\sum_\bj \sum_{\bx_\bk\in \ST(\bx_\bj)}u(\bx_\bk)\,
\Pol_{\bj,\bk}\Big(\frac{\bx-\bx_{\bj}}{h_{\bj}}\Big)\,
{\ee}^{- |\bx-\bx_{\bj}|^{2}/h_{\bj}^2} \, .
\end{equation}
Using  the following lemma it is easy to derive cubature formulae 
for \eqref{defint}.
\begin{lemma} \label{ptos}
For any
$ \ds
\Pol(\bx)=\sum_{|\bb|=0}^{L}    c_{\bb}\, \bx^\bb
$
one can write
$ \ds
\Pol(\bx) \ee^{-|\bx|^2}
=\sum_{|\bb|=0}^{L} c_{\bb} \,
\Ss_{\bb}(\partial_\bx) \, \ee^{-|\bx|^2}
$
with the polynomial $\Ss_{\bb}(\btt)$ being defined by
\begin{equation}\label{e:S}
\Ss_{\bb}(\btt)=\Big(\frac{1}{2\, i}\Big)^{|\bb|} H_{\bb}
\Big(\frac{\btt}{2\, i} \Big)\, ,
\end{equation}
where $H_{\bb}$ denotes the Hermite polynomial of $n$ variables
$
H_\bb(\btt)={\ee}^{\,|\btt|^2}( - \partial_\btt )^{\bb}
{\ee}^{-|\btt|^2}
$.
\end{lemma}
\textit{Proof.} We are looking for the polynomial $\Ss_{\bb}(\btt)$
defined by the relation
\begin{equation}\label{delta}
\bx^{\bb}
\ee^{-|\bx|^{2}}=\Ss_{\bb}(\partial_\bx) \ee^{-|\bx|^{2}},\quad \bx \in \R^{n}.
\end{equation}
Since
\[
\Ff(\Ss_{\bb}(\partial_\bx)
\ee^{-|\bx|^{2}})(\l)=\pi^{n/2} \ee^{-\pi^{2}|\l|^{2}}
\Ss_{\bb}(2\pi i \l)
\]
and
\[
\Ff(\bx^{\bb}
\ee^{-|\bx|^{2}})(\l)=\pi^{n/2} \Big(-\frac{\partial_\l}{2 \pi
i}  \Big )^{\bb} \ee^{-\pi^{2} |\l|^{2}}
\]
we obtain (\ref{e:S}).\hfill  $\blacksquare$
\bigskip

In view of Lemma \ref{ptos} we can write
$
\Pol_{\bj,\bk}(\bx)\, {\ee}^{-|\bx|^{2}}=
\T_{\bj,\bk} (\partial_\bx)\, {\ee}^{-|\bx|^{2}} 
$
with some polynomials $\T_{\bj,\bk}(\bx)$.
Then
(\ref{mainM2}) can be rewritten as
\[
M u(\bx)= \sum_\bj \sum_{\bx_\bk\in \ST(\bx_\bj)}u(\bx_\bk)\,
\T_{\bj,\bk} (-h_{\bj}\, \partial_{\bx_\bj})\,
{\ee}^{-|\bx-\bx_{\bj}|^2/h_{\bj}^{2}} \, .
\]
The cubature formula for the integral  $\K u$ is obtained by
replacing $u$ by its quasi-interpolant $M u$
\begin{equation} \label{cubat}
\tilde{\K}u(\bx)=\K M u(\bx)=
\sum_\bj \sum_{\bx_\bk\in \ST(\bx_\bj)}u(\bx_\bk)\,
\T_{\bj,\bk}(-h_{\bj} \, \partial_{\bx_\bj})
\,h_{\bj}^{n}\,\int_{\R^{n}} g(h_{\bj}|\bz|)\,
{\ee}^{-|\bz+\btt_{\bj}|^{2}}\, d \bz \, ,
\end{equation}
where $\btt_{\bj}=(\bx-\bx_{\bj})/h_{\bj}$.
By introducing
spherical coordinates in $\R^{n}$ we obtain
\[
\int_{\R^{n}}\!\! g(h_{\bj}|\bz|)  \ee^{-|\bz+\btt_{\bj}|^{2}}d
\bz=
\ee^{-|\btt_{\bj}|^{2}} \int_{0}^{\infty} \r^{n-1} g(h_{\bj}
\r) \, \ee^{-\r^{2}} d\,  \r \int_{S^{n-1}} \ee^{-2\r |\btt_{\bj}|
\cos(\o_{\btt_{\bj}},\o)} \, d \s_{\o} \, ,
\]
where $S^{n-1}$ is the unit sphere in $\R^{n}$.
The integral
over $S^{n-1}$ can be represented by means of the
modified Bessel functions of the first kind $I_{n}$ in the following way
\[
\int_{S^{n-1}} \ee^{-2\r |\btt_{\bj}|\cos(\o_{\btt_{\bj}},\o)}d \s_{\o}=
\frac{2\, \pi^{(n-1)/2}}{\G(\frac{n-1}{2})}
\int_{0}^{\pi} \ee^{-2 \r |\btt_{\bj}| \cos \th } (\sin \th)^{n-2} \, d \th=
2 \pi^{n/2} (\r \,
|\btt_{\bj}|)^{1-n/2} I_{\frac{n-2}{2}}(2\r |\btt_{\bj}|)
\]
(see \cite[p.154]{SW} and \cite[p.79]{Wa}).
If we denote by
\[
\Ll(r)=2\, \pi^{n/2}\,
r^{1-n/2}\ee^{-r^{2}}
\int_{0}^{\infty} \r^{n/2} \, \ee^{-\r^{2}} g(h_{\bj} \, \r)
\, I_{(n-2)/2}(2\r \, r) \, d\r \, ,
\]
then \eqref{cubat} leads to  
the following cubature formula for
the integral $Ku$ 
\[
\tilde{\K} u(\bx)=\sum_\bj h_{\bj}^{n} \sum_{\bx_\bk\in
\ST(\bx_\bj)}u(\bx_\bk)\,
\T_{\bj,\bk} (-h_{\bj} \, \partial_{\bx_\bj})
\,\Ll\Big( \frac{|\bx-\bx_{\bj}|}{h_{\bj}}\Big)
\, . 
\]

\paragrafo{Construction of the $\Th$-function with Gaussians}\label{S:QUA}
\setcounter{equation}{0}
In this section we propose 
a method to construct the approximate
partition of unity for the basis functions
\[
\hj(\bx)=(\pi\,D)^{-n/2} \ee^{-|{\bx}-\bx_\bj|^{2}/h_\bj^2 D} 
\]
if the set of nodes $\{\bx_\bj\}$ satisfy
Condition \ref{cond1} piecewise with different grid sizes $h_j$.

\subsection{Scattered nodes close to a piecewise uniform grid}\label{SS:61}
Let us explain the assumption on the nodes: 
Suppose that 
a subset of nodes $\bx_{\bj} \in J_0$ 
satisfies Condition \ref{cond1} with $h=h_1$.
The remaining nodes $\bx_\bk \in  \bX \setminus J_0$ 
lie in a bounded domain $\Omega \subset \R^n$ and satisfy 
Condition \ref{cond1} with $h=h_2=H h_1$  for some small $H$. 
To keep good local properties of quasi-interpolants
one wants to approximate the data at these nodes by
functions of the form polynomial times
$\ee^{-|{\bx}-\bx_\bk|^{2}/h_2^2 D}$, whereas outside $\Omega$
quasi-interpolants with functions of the form polynomial times
$\ee^{-|{\bx}-\bx_\bj|^{2}/h_1^2 D}$ should be used.

Our aim  is, to develop
a simple method 
to construct 
polynomials $\Pol_\bj$
such that
\begin{equation}\label{th0}
\Th(\bx)=
(\pi\, D)^{-n/2} \Big( \sum_{\bx_\bj \in J_1} 
\Pol_\bj \Big( \frac{\bx-\bx_\bj}{h_1 \dd} \Big)
\ee^{-|{\bx}-\bx_\bj|^{2}/h_1^2 D} +
\sum_{\bx_\bk \in J_2} 
\Pol_\bk \Big( \frac{\bx-\bx_\bk}{h_2 \dd} \Big)
\ee^{-|{\bx}-\bx_\bk|^{2}/h_2^2 D}\Big)
\end{equation}
is almost the constant function $1$.
Here $J_2$ denotes the set of  nodes $\bx_\bk \in \Omega$
and $J_1 = \{\bx_\bj\} \setminus J_2$ the remaining nodes. 

First we derive a piecewise uniform grid on $\R^n$ which is
associated to the splitting of the set of scattered nodes into $J_1$ and $J_2$.
We start with
Poisson's summation formula for Gaussians
\[
(\pi \, D)^{-n/2}\sum_{\bm \in \Z^n}
\ee^{-|{\bx}-h_1\bm|^{2}/h_1^2D} 
= \sum_{\bk \in \Z^n} \ee^{-\pi^2 D |\bk|^2} \ee^{\, 2 \pi i
  (\bx,\bk)/h_1} \, ,
\]
which shows that
\[
\Big|1- (\pi \, D)^{-n/2} \sum_{\bm \in \Z^n}
\ee^{-|{\bx}- h_1\bm|^{2}/h_1^2 D}
\Big| \le C_1 \, \ee^{-\pi^2 D}
\]
with some constant $C_1$ depending only on the space dimension.

Thus for any $\e>0$ there exists $D>0$ such that
the function system $\{\ee^{-|{\bx}-h_1\bm|^{2}/h_1^2 D}\}_{\bm \in \Z^n}$ 
forms an approximate partition of unity
with accuracy $\e$. 
We can represent any of these functions very accurately by a linear combination of
dilated Gaussians due to the equation (see \cite{MS5})
    \begin{equation} \label{refin}
        \begin{split}
\ee^{- |\bx|^2/{D_1}} =
\Big( &\frac{{D}_1} { \pi {D}({D}_1-h^2 {D})} \Big)^{n/2}
    \sum_{\bm \in {\Z}^n}  
        \ee^{\,- {h^2 |\bm|^2}/({D}_1-h^2 {D})}
    \, \ee^{\,- |\bx-h \bm|^2/h^2 D } \\
  &- \ee^{- |\bx|^2/{D_1}} 
        \sum_{\bk \in {\Z}^n \backslash \{\boldsymbol{0}\}} 
        \ee^{\, 2 \pi i ({D}_1-h^2 D ) (\mathbf{x}, \bk ) /h {D}_1 } \,
        \ee^{\, - \pi^2 {D} ({D}_1-h^2 D)
|\bk|^2/{{D}_1} } ,
        \end{split}
    \end{equation}
which is valid for any ${D}_1 > h^2 {D} > 0$.
Applied to our setting  with $h=h_2$ and  ${D}_1= h_1^2 D$
we obtain the approximate refinement relation
    \begin{equation} \label{refinest}
\Big|
\ee^{- |\bx|^2/ h_1^2 D} -
    \sum_{\bk \in {\Z}^n}
a_\bk     \, \ee^{- |\bx -h_2 \bk |^2/{ h_2^2 D}} \Big| 
\le C_2 \, \ee^{-  |\bx|^2/h_1^2 D} \, \ee^{- \pi^2 D(1-H^2)} 
    \end{equation}
(because by assumption $h_2 = H h_1$) with the coefficients
\[
a_\bk   =  \big({ \pi  D(1- H^2)} \big)^{-n/2}
\ee^{- H^2 |\bk|^2/ (1-H^2)D} \, .
\]
Again, the constant $C_2$ depends only on the space dimension.
Define by $S \in \Z^n$ the minimal index set such that
\[
\sum_{\bk \in {\Z}^n \setminus S}
a_\bk < \ee^{- \pi^2 D (1-H^2)} \, .
\]
Then  it is clear from \eqref{refinest} that
for any disjoint  $Z_1$ and $Z_2$ with
$Z_1 \cup Z_2 = \Z^n$ 
    \begin{equation} \label{twoscale} 
        \begin{split}
\Big|1- (\pi \, D)^{-n/2} \Big( \sum_{\bm \in Z_1}
\ee^{-|{\bx}- h_1\bm|^{2}/h_1^2 D} 
+ \sum_{\bm \in Z_2} \sum_{\bk \in S} 
a_\bk  & \,
\ee^{- |\bx -h_1 \bm -h_2 \bk |^2/{ h_2^2 D}} \Big) \Big| 
\le C_3 \, \ee^{- \pi^2 D(1-H^2)} \, .
        \end{split}
    \end{equation}

\begin{cond} \label{cond4}
{\em 
Denote $Z_2 = \{\bm \in \Z^n: h_1\bm + h_2 \bk \in
\Omega \; \textrm{for all} \; \bk \in S\}$. 
The constant $\k_1$ of Condition~\ref{cond1} and the domain $\Omega$
are such that for all
nodes $\bx_\bk \in \Omega$, i.e. the nodes belonging to $J_2$,
one can find $\bm \in Z_2$, $\bk \in S$ with
 $|\bx_\bk - h_1\bm - h_2 \bk|< \kappa_1 h_2$.}
\end{cond}
Setting $Z_1=\Z^n \setminus Z_2$ we connect  the index sets 
$Z_1$, $Z_2$ 
with  the splitting of the scattered nodes into $J_1$, $J_2$.
By this way we construct an approximate partition of unity
using Gaussians
with the "large" scaling factor $h_1$ 
centered  at the uniform grid $G_1:=\{h_1 \bm\}_{\bm \in Z_1}$
outside
$\Omega$ and
using Gaussians with scaling factor $h_2$ and the centers 
$G_2:= \{h_1\bm + h_2 \bk\}_{\bm \in Z_2, \bk \in S}$ in
$\Omega$. 

It is obvious, that the above definition of piecewise quasi-uniformly
distributed scattered nodes and the construction of an associated
approximate partition of unity on piecewise uniform grids
can be extended to finitely  many scaling factors $h_j$.
Since there will be no difference for the subsequent considerations 
we will restrict to the two-scale case.

From \eqref{twoscale} we see that for
any $\e > 0$, and given $ h_1$ and $h_2$ there exists $D >0$ such that
the linear combination 
    \begin{equation} \label{twoscale1} 
(\pi \, D)^{-n/2} \Big( \sum_{g_1 \in G_1}
\ee^{-|{\bx}- g_1|^{2}/h_1^2 D} 
+ \sum_{g_2 \in G_2 }
\tilde a_{g_2}  \,
\ee^{- |\bx - g_2|^2/{ h_2^2 D}}\Big)
    \end{equation}
with $\tilde a_{g_2} = a_\bk$ for 
$g_2=h_1\bm + h_2 \bk, \,\bm \in Z_2, \bk \in S$, 
approximates the constant function $1$ with an error less than $\e/2$.
The  
idea of constructing the $\Th$-function \eqref{th0}
is to choose  for each  $g_1 \in G_1$ and $g_2 \in G_2$
finite sets of nodes
$\Sigma(g_1) \subset J_1$
and $\Sigma(g_2) \subset J_2$, respectively,
and to determine 
polynomials $\Pol_{\bj,g_\ell}$
such that
\[
\sum_{\bx_\bj \in  \Sigma(g_\ell)}\Pol_{\bj,g_\ell}
\Big( \frac{\bx-\bx_\bj}{h_\ell\dd} \Big)
\ee^{-|{\bx}-\bx_\bj|^{2}/h_\ell^2 D} \quad \mbox{approximate} \quad 
\ee^{-|{\bx}-g_\ell|^{2}/h_\ell^2 D} \, , \; \ell =1,2 \, .
\]
If the $L_\infty$-error of the sums over $g_\ell$
can be controlled,
then we get
\begin{equation*}
\sum_{g_1 \in G_1}
\ee^{-|{\bx}-g_1|^{2}/h_1^2 D} 
\asymp\sum_{\bx_\bj \in J_1}  
\Pol_\bj\Big( \frac{\bx-\bx_\bj}{h_1 \dd} \Big) 
\ee^{-|{\bx}-\bx_\bj|^{2}/h_1^2 D} 
\end{equation*}
with the polynomials
\begin{equation} \label{polm}
\Pol_\bj  = \sum_{g_1 \in G(\bx_\bj)}
\Pol_{\bj,g_1} 
\end{equation}
and 
\begin{equation*}
\sum_{g_2 \in G_2} \tilde a_{g_2}
\ee^{-|{\bx}-g_2|^{2}/h_2^2 D} 
\asymp\sum_{\bx_\bk \in J_2}  \Pol_\bk\Big( \frac{\bx-\bx_\bk}{h_2 \dd} \Big) 
\ee^{-|{\bx}-\bx_\bk|^{2}/h_2^2 D} 
\end{equation*}
with the polynomials
\begin{equation} \label{polm1}
\Pol_\bk  = \sum_{g_2 \in G(\bx_\bk)}
\tilde a_{g_2} \Pol_{\bk,g_2} \, ,
\end{equation}
where we denote $G(\bx_\bj) = \{g :
\bx_\bj \in \Sigma(g)\}$.
Note that we have to choose the subsets $\Sigma(g_\ell)$ such that
the sets $G(\bx_\bj) \subset G_\ell$ are finite and nonempty
for any node $\bx_\bj \in J_\ell$.
Additionally, one has to choose these sets such that
for some $\kappa_1 > 0$ and any 
$g_\ell \in G_\ell$ the ball $B(g_\ell, \kappa_1 h_\ell)$
contains at least one node $\bx_\bj \in J_\ell$.
This is always possible, 
since  Conditions~\ref{cond1} resp. \ref{cond4}
are valid.

The proposed construction method of $\Th$ does not require
solving a large algebraic system. Instead, to obtain
the local representation of $\Th$ one has to solve 
a small number of approximation problems, which are reduced in the
next sections to linear systems of moderate size.

After this preparation we 
write $\Th$ as
    \begin{equation*}
        \begin{split}
\Th(\bx)= &(\pi \, D)^{-n/2} \sum_{g_1 \in G_1}
\ee^{-|{\bx}- g_1|^{2}/h_1^2 D} 
+ \sum_{g_1 \in G_1} \o_{g_1} (\frac{\bx}{h_1}) \\
&+ (\pi \, D)^{-n/2} \sum_{g_2 \in G_2 }
\tilde a_{g_2}  \,
\ee^{- |\bx - g_2|^2/ h_2^2 D} 
+ \sum_{g_2 \in G_2} \tilde a_{g_2} \o_{g_2} (\frac{\bx}{h_2}) \, ,
        \end{split}
    \end{equation*}
where
\begin{equation}\label{omegam}
\o_{g_\ell} (\by)=(\pi \, D)^{-n/2}
\Big\{
\sum_{h_\ell \by_\bj  \in \Sigma(g_\ell)}
\Pol_{\bj,g_\ell}\Big( \frac{\by-\by_\bj}{\dd} \Big)
\ee^{-|\by-\by_\bj|^{2}/D} 
- \ee^{-|{\by}- g_\ell/h_\ell|^{2}/D} 
\Big\}
\end{equation}
with $\by_{\bj} = \bx_{\bj}/h_\ell$, $\bx_{\bj} \in J_\ell$. 
Hence for sufficiently large $D$
\begin{equation}\label{thm1}
\big|\Th(\bx)-1\big| <\frac{\e}{2} + \sum_{g_1 \in G_1}
\Big|\o_{g_1} \big(\frac{\bx}{h_1}\big)\Big| +
\sum_{g_2 \in G_2} \tilde a_{g_2} 
\Big|\o_{g_2} \big(\frac{\bx}{h_2}\big)\Big| \, .
\end{equation}

\subparagrafo{Construction of Polynomials}
Let us introduce 
\begin{equation} \label{omegam_2}
\o(\by):=(\pi \, D)^{-n/2}
\Big\{
\sum_{\by_\bj \in \Sigma}
\Pol_{\bj}\Big( \frac{\by-\by_\bj}{\dd} \Big)
\ee^{-|\by-\by_\bj|^{2}/D} 
- \ee^{-|{\by}|^{2}/D} 
\Big\} ,
\end{equation}
where $\Sigma$ is some finite point set in $\R^n$. 
We will describe a method
for constructing 
polynomials $\Pol_{\bj}$ such that
$\ee^{\, \rho|\by|^2}|\o(\by)|$ for some $\rho>0$ becomes small.
In what follows we use the representation
\[
\Pol_{\bj}(\bx)=\sum_{|\bb|=0}^{L_\bj} c_{\bj,\bb}\, \bx^\bb \, .
\]
Hence by Lemma \ref{ptos}
\[
\Pol_{\bj}\Big( \frac{\by-\by_\bj}{\dd} \Big)
\ee^{-|\by-\by_\bj|^{2}/D} = 
\sum_{|\bb|=0}^{L_{\bj}} c_{\bj,\bb} \Ss_{\bb}(\dd \partial_\by) \,
\ee^{- |\by-\by_{\bj}|^2/D} \, .
\]
and $\o$
can be written as
\begin{equation}\label{4.1}
\o (\by) = 
 {(\pi D)^{-n/2}}
\Big(\sum_{ \by_\bj \in \Sigma} \sum_{|\bb|=0}^{L_{\bj}} c_{\bj, \bb}  
\Ss_{\bb}(\dd \partial_\by) \,
\ee^{- |\by-\by_{\bj}|^2/D}-
\ee^{-|\by|^2/D}\Big).
\end{equation}
To estimate the $L_\infty$-norm  of $\o$ we
represent this function as  convolution.
\begin{lemma}\label{DUE}
Let $\Pol(\btt)$ be a polynomial and let   $0<D_0<D$. Then
\[
\Pol(\partial_\bxi) \, \ee^{-{|\bxi-\bx|^{2}}/{D}}=c_{1}\,
\ee^{-{|\bxi|^{2}}/{(D-D_{0})}}\ast
\Pol(\partial_\bxi)\,\ee^{-{|\bxi-\bx|^{2}}/{D_{0}}} \, ,
\]
where $\ast$ stands for the convolution operator and
\[
c_{1}=\Big(\frac{D}{\pi D_{0}(D-D_{0})}\Big)^{n/2} \, .
\]
\end{lemma}
\textit{Proof.} From
\[
\ee^{-{|\bxi-\bx|^{2}}/{D}}=c_1\,
\int_{\R^{n}} \ee^{-{|\bxi-\btt|^{2}}/{(D-D_{0})}} \>
\ee^{-{|\btt-\bx|^{2}}/{D_{0}}} \, d \btt  
\]
we obtain
\[
\Pol(\partial_\bxi)\ee^{- {|\bxi-\bx|^{2}}/{D}}=
\Pol(-\partial_\bx)\ee^{-{|\bxi-\bx|^{2}}/{D}}=
c_1
\int_{\R^{n}} \ee^{-{|\bxi-\btt|^{2}}/{(D-D_{0})}}
\Pol(\partial_\btt)\ee^{-{|\btt-\bx|^{2}}/{D_{0}}} \, d
\btt \, .\qquad \blacksquare
\]
\medskip

Using Lemma \ref{DUE} and 
\eqref{4.1}  we  write 
$\o$
as 
\begin{equation}\label{ast}
\o(\by)=
\frac{c_{1}}{(\pi
D)^{n/2}}
\int_{\R^n}
\ee^{- |\by-\btt|^{2}/(D-D_{0})}
\Big( \sum_{\by_\bj \in \Sigma} \sum_{|\bb|=0}^{L_\bj}
c_{\bj,\bb} \Ss_{\bb}(\dd \partial_\btt)
\ee^{-|\btt-\by_{\bj}|^{2}/D_{0}}-
\ee^{-|\btt|^{2}/D_{0}} \!\Big)  \,d \btt
\end{equation}
and, by Cauchy's inequality, we obtain 
\begin{equation} \label{simpast}
||\o||_{L_\infty}\leq
c_{2}\, \Big\Vert \sum_{\by_\bj \in \Sigma} \sum_{|\bb|=0}^{L_\bj}
c_{\bj,\bb}  \Ss_{\bb}(\dd \partial_\btt)
\ee^{-|\btt-\by_{\bj}|^{2}/D_{0}}-
\ee^{-|\btt|^{2}/D_{0}}\Big\Vert_{L^{2}} \, ,
\end{equation}
where 
\[
c_{2}=(\pi D_{0})^{-n/2}(2 \pi (D-D_{0}))^{-n/4}\, .
\]
If we define polynomials $T_{\bb}$  by
\begin{equation} \label{defT}
T_{\bb}( \bx)
 = \ee^{\,|\bx|^{2}/D_0} \Ss_{\bb}(\dd \partial_\bx)
\ee^{-|\bx|^{2}/D_0} \, ,
\end{equation}
then
\begin{equation*}
||\o||_{L_\infty}\leq
c_{2}\, \Big\Vert\ee^{-|\,\cdot \,|^{2}/D_{0}}-
 \sum_{\by_\bj \in \Sigma} \sum_{|\bb|=0}^{L_\bj}
c_{\bj,\bb}  T_{\bb}(\,\cdot \,-\by_{\bj}) 
\ee^{-|\,\cdot \, -\by_{\bj}|^{2}/D_{0}} \Big\Vert_{L^{2}} \, .
\end{equation*}
An estimate for  the sum of $|\o_{g_\ell}|$
can be derived from
\begin{lemma} \label{constr}
Let $0<D_0<D$ and denote $\ds  \rho = \frac{D-D_0}{(D-D_0)^2+DD_0}$.
Then the estimate
\begin{equation}  \label{ast1}
\begin{split}
\sup_{\R^n} & \, |\o (\by)|\, \ee^{\, \rho|\by|^2} \,
\leq c_3 \sqrt{Q(\cc)}
\end{split}
\end{equation}
is valid, where for 
$\cc =\{c_{\bj,\bb}\}$
the quadratic form $Q(\cc)$ is defined by
\begin{equation} \label{defQ}
Q(\cc)  = \!
\int_{\R^n}
\ee^{2(D-D_0)|\btt|^2/D D_0}
\Big(\ee^{-|\btt|^2/D_0}- \! \! \sum_{\by_\bj \in \Sigma} 
\sum_{|\bb|=0}^{L_{\bj}} 
c_{\bj,\bb} T_{\bb}(\btt-\by_\bj)
\,\ee^{-|\btt-\by_\bj|^2/D_0} \Big)^2  d \btt 
\end{equation}
and
\[
c_3 =
\frac{D^{n/4} }
{(2 \pi^3 D_0 (D-D_0)((D-D_0)^2+D D_0))^{n/4}} \, .
\]
\end{lemma}
{\bf Proof.} 
Starting with \eqref{ast}, using \eqref{defT} and
\[
|\bx-\btt|^2 = \Big| \sqrt{a}\bx  - \frac{\btt}{ \sqrt{a}} \Big|^2+
(1-a)|\bx|^2 + \frac{a-1}{a} |\btt|^2
\]
for $a>0$, we derive the representation
\begin{equation*}
\begin{split}
\o(\by)  = \frac{c_{1}}{(\pi
D)^{n/2}}\, \ee^{-(1-a)|\by|^2/(D-D_0)} 
&\int_{\R^n}  
\ee^{-|\btt- a \bx|^2/a(D-D_0)} \, 
\ee^{(1-a)|\btt|^2/a(D-D_0)} \\
&\times 
\Big(\sum_{\by_\bj \in \Sigma}\sum_{|\bb|=0}^{L_{\bj}} 
c_{\bj,\bb}
T_{\bb}(\btt-\by_\bj)
\ee^{-|\btt-\by_\bj|^{2}/D_0} -\ee^{-|\btt|^2/D_0} ) \, d \btt \, .
\end{split}
\end{equation*}
Then Cauchy's inequality leads to 
\begin{equation} \label{wo1} 
\begin{split}
\Big|&\,\o(\by)\, \ee^{(1-a)|\by|^2/(D-D_0)}\Big| \\
&\le
c_3 \Big(\int_{\R^n}
\ee^{ 2(1-a)|\btt|^2/a(D-D_0)}
\Big(\sum_{ \by_\bj \in \Sigma}\sum_{|\bb|=0}^{L_{\bj}} 
c_{\bj,\bb}
T_{\bb}(\btt-\by_\bj)
e^{-|\btt-\by_\bj|^{2}/D_0} -\ee^{-|\btt|^2/D_0}\Big)^2 d \btt \Big)^{1/2}
\end{split}
\end{equation}
with
\begin{equation*} 
\begin{split}
c_3 &= 
 (\pi D_0)^{-n/2}
\Big(\frac{a }{2 \pi (D-D_0)}\Big)^{n/4} \, .
\end{split}
\end{equation*}
If we choose the parameter $a$ such that
\[
\frac{(1-a)|\btt|^2}{a(D-D_0)} - \frac{|\btt|^2}{D_0} =
-\frac{|\btt|^2}{D \,} \, ,
\quad 
\mbox{i.e.} \quad
a = \frac{D\,D_0}{(D-D_0)^2+D\,D_0} \, ,  
\]
then the right hand side of \eqref{wo1} takes the form \eqref{defQ}. 
\hfill $\blacksquare$
\bigskip

Next we estimate 
\begin{equation} \label{defrm}
r := \min_{\cc} Q (\cc) \, .
\end{equation}
Using \eqref{defT}, after elementary calculations one  obtains 
\begin{equation*}
\begin{split}
&Q(\cc)=
\int_{\R^n}
\Big( \ee^{(D-D_0)|\btt|^2/D D_0}
\Big(\sum_{\by_\bj \in \Sigma} \sum_{|\bb|=0}^{L_{\bj}} 
c_{\bj,\bb} 
\Ss_{\bb}(\dd \partial_\btt) \,\ee^{-|\btt-\by_\bj|^2/D_0}
- \ee^{-|\btt|^2/D_0}\Big) \Big)^2 \, d \btt \\
&= \Big(\frac{\pi D}{2}\Big)^{n/2}\Big( 1 -2
\sum_{\by_\bj \in \Sigma} \sum_{|\bb|=0}^{L_{\bj}} 
c_{\bj,\bb} \C_{\bb,0}(\by_\bj,\boldsymbol{0})
+\sum_{\by_\bj ,\by_\bk\in \Sigma} \sum_{|\bb|=0}^{L_{\bj}} 
\sum_{|\bg|=0}^{L_{\bk}}c_{\bj,\bb}c_{\bk,\bg}
\C_{\bb,\bg}(\by_\bj,\by_\bk)\Big) 
\end{split}
\end{equation*}
with 
\[
\C_{\bb,\bg}(\bx,\by) := 
\Ss_{\bb}(- \dd \partial_{\bx})
\Ss_{\bg}(-\dd \partial_{\by}) \,
\ee^{\,(D-D_0)(|\bx|^2+|\by|^2)/D^2_0} \, 
\ee^{ -D|\bx-\by|^2/2 D^2_0} \, .
\]
The minimum of $Q(\cc)$ is attained by the solution 
$\cc=\{c_{\bj,\bb}\}$
of the linear system
\begin{equation}\label{sys1}
\sum_{\by_\bj \in \Sigma} \sum_{|\bb|=0}^{L_{\bj}} 
\C_{\bb,\bg}(\by_{\bj},\by_{\bk}) c_{\bj,\bb}=
\C_{\boldsymbol{0},\bg}(\boldsymbol{0},\by_{\bk})\,, 
\quad  \by_{\bk} \in \Sigma \, ,\;
0 \le |\bg| \le L_{\bk} \, .
\end{equation}
Then by Lemma \ref{constr} the sum
\begin{equation} \label{lsdef}
\sum_{\by_\bj \in \Sigma}\Pol_{\bj}\Big( \frac{\by-\by_\bj}{ \dd} \Big)
\ee^{-|{\by}-\by_\bj|^{2}/ D}=\sum_{\by_\bj \in \Sigma} 
\sum_{|\bb|=0}^{L_{\bj}}
c_{\bj,\bb} \Big( \frac{\by-\by_\bj}{ \dd} \Big)^\bb
\ee^{-|\by -\by_{\bj}|^{2}/D}
\end{equation}
approximates  $\ee^{-|\by|^2/D}$ 
with
\[
(\pi \, D)^{-n/2}\Big|\ee^{-|\by|^2/D}-\sum_{\by_\bj \in
  \Sigma}\Pol_{\bj}
\Big( \frac{\by-\by_\bj}{ \dd} \Big)
\ee^{-|{\by}-\by_\bj|^{2}/ D}\Big|\le c_3 \, \ee^{\,- \rho|\by|^2}
r^{1/2} \, .
\]
In the next section we show that \eqref{sys1}
has a unique solution and give an estimate of $r$.

\subparagrafo{Existence and estimates} 
Let us  give another representation of
the quadratic form
$Q (\cc)$ defined by \eqref{defQ}.
Introduce
the transformed points
\[
\btt_{\bj} = \frac{D}{D_0} \by_{\bj}  \, , \quad
\by_{\bj} \in \Sigma \, ,
\]
then, because of 
\begin{equation*}
\begin{split}
\frac{D-D_0}{D\, D_0} |\btt|^2 - \frac{1}{D_0}|\btt - \by_{\bj}|^2
= -\frac{1}{D}|\btt-\btt_{\bj}|^2 
+ \frac{D-D_0}{D_0^2} \, |\by_{\bj}|^2
\end{split}
\end{equation*}
$Q(\cc)$ can be written as
 \begin{equation} \label{ast1n}
 \begin{split}
Q(\cc)= \int_{\R^n}\Big( \ee^{-|\btt|^2/D} - 
\! \! \sum_{\by_\bj \in \Sigma} 
 \ee^{(D-D_0)|\by_{\bj}|^2/D_0^2} 
 \sum_{|\bb|=0}^{L_{\bj}} c_{\bj,\bb} 
 T_{\bb}(\btt -\by_{\bj}) \,\ee^{-|\btt-\btt_{\bj}|^2/D}
\Big)^2 d \btt \, .
\end{split}
\end{equation}
Since $T_{\bb}$ are polynomials of degree $\bb$, the minimum problem for $Q(\cc)$
is equivalent
to finding the best $L_2$-approximation
\[
\min_{d_{\bj,\bb}} \; \int_{\R^n}\Big( \ee^{-|\btt|^2/D} - 
 \sum_{\by_\bj \in \Sigma} 
 \sum_{|\bb|=0}^{L_{\bj}} d_{\bj,\bb} 
 (\btt -\btt_{\bj})^\bb \,\ee^{-|\btt-\btt_{\bj}|^2/D}
\Big)^2 d \btt \, .
\]
\begin{lemma} \label{lsexist}
Let $\{\bx_\bj\}$ a finite collection of nodes.
For all $L_{\bj} \ge 0$ the  polynomials
$\Pol_{\bj}$ of degree $L_{\bj}$, which minimize
\begin{equation} \label{defnorm}
\Big\| \ee^{-|\, \cdot\,|^2} - 
\sum_{\bj} \Pol_{\bj}(\,\cdot\,-\bx_\bj) \, \ee^{-|\,\cdot\,-\bx_\bj|^2} 
\Big\|_{L_2} \, ,
\end{equation}
are uniquely determined.
\end{lemma}
\textit{Proof.} 
The application of Lemma \ref{ptos} gives for 
$
\ds  \Pol_{\bj}(\bx)=\sum_{|\bb|=0}^{L_{\bj}}    c_{\bj,\bb}\, \bx^\bb
$
\begin{equation*}
\begin{split}
&\Big\| \ee^{-|\, \cdot\,|^2} - 
\sum_{\bj} \Pol_{\bj}(\,\cdot\,-\bx_\bj) \ee^{-|\,\cdot\,-\bx_\bj|^2} 
\Big\|^2_{L_2}
= \int_{\R^n} \Big(\ee^{-|\bx|^2} - 
\sum_{\bj}  \sum_{|\bb|=0}^{L_{\bj}} c_{\bj,\bb} \,
\Ss_{\bb}(\partial_\bx) \,\ee^{-|\bx-\bx_\bj|^2}\Big)^2 d \bx
\\
&=\Big(\frac{\pi}{2}\Big)^{n/2}\Big( 1 - 2\sum_{\bj} 
\sum_{|\bb|=0}^{L_{\bj}} c_{\bj,\bb} \B_{\bb,0}(\bx_{\bj},0)
+\sum_{\bj,\bk} \sum_{|\bb|,|\bg|=0}^{L_{\bj},L_{\bk}} c_{\bj,\bb}c_{\bk,\bg}
\B_{\bb,\bg}(\bx_{\bj},\bx_{\bk})\Big)  ,
\end{split}
\end{equation*} 
where we use the notation
\[
\B_{\bb,\bg}(\bx,\by)=
   \Ss_{\bb}(-  \partial_\bx) \,
\Ss_{\bg}(- \partial_\by)\, \ee^{- |\bx -\by|^2/2} \, .
\]
The coefficients $\{\cc_{\bj,\bb}\}$ are chosen  to minimize \eqref{defnorm},
that is the vector
$\{\cc_{\bj,\bb}\}$
is a solution of the linear system
\begin{equation}\label{sys}
\sum_{\bj} \sum_{|\bb|=0}^{L_{\bj}} 
 c_{\bj,\bb} \, \B_{\bb,\bg}(\bx_{\bj},\bx_{\bk})=
\B_{0,\bg}(0,\bx_{\bk}).
\end{equation}
To show that the matrix of this system is positive definite
we use the representation
\[
\ee^{- |\bx -\by|^2/2}=
({2\pi})^{-n/2}\int_{\R^n}
\ee^{-|\btt|^2/2} \ee^{i(\btt,\bx)}  \ee^{- i(\btt,\by)} \, d \btt  \, ,
\]
which implies
\[
\B_{\bb,\bg}(\bx,\by)=
({2\pi})^{-n/2}\int_{\R^n}
\Ss_{\bb} ( - i \btt) \, 
\ov{\Ss_{\bg} ( - i \btt)} \, 
\ee^{-|\btt|^2/2} \ee^{i(\btt,\bx)}  \ee^{- i(\btt,\by)} \, d \btt  \, .
\]
Let $\{v_{\bj,\bb}\}$ be a constant vector and consider
the sesquilinear form
\begin{equation*}
\begin{split}
\sum_{\bj,\bk} & \sum_{|\bb|,|\bg|=0}^{L_{\bj},L_{\bk}}
\B_{\bb,\bg}(\bx_{\bj},\bx_{\bk}) \,v_{\bj,\bb} \, \ov {v_{\bk,\bg}} \\
&= ({2\pi})^{-n/2}\sum_{\bj,\bk} 
\sum_{|\bb|,|\bg|=0}^{L_{\bj},L_{\bk}}
v_{\bj,\bb} \, \ov{v_{\bk,\bg}} \int_{\R^n}
\Ss_{\bb} ( - i \btt) \,
\ov{\Ss_{\bg} ( - i \btt)}\, 
\ee^{-|\btt|^2/2} \ee^{ \, i(\btt,\bx_\bj-\bx_\bk)}
\, d \btt \\
&= ({2\pi})^{-n/2} \int_{\R^n}
\ee^{-|\btt|^2/2} \Big|\sum_{\bj} \sum_{|\bb|=0}^{L_\bj}
v_{\bj,\bb} \, \Ss_{\bb} ( - i \btt) \, \ee^{i(\btt,\bx_\bj)}\Big|^2
\, d \btt \ge 0 \, .
\end{split}
\end{equation*} 
The change of integration and summation is valid because the
integrand is absolutely integrable and the sums are finite.
We have to show that the inequality is strict when $\{v_{\bj,\bb}\} \neq 0$.
This is equivalent to show that
\begin{equation*}
\begin{split}
\sigma(\btt) &=\sum_{\bj} \sum_{|\bb|=0}^{L_\bj}
v_{\bj,\bb} \, \Ss_{\bb} ( - i \btt) \, \ee^{i(\btt,\bx_\bj)} = 0
\end{split}
\end{equation*} 
identically only if all components
$v_{\bj,\bb} = 0$ for all $\bj$ and $\bb$.
To this end similar to \cite[Lemma~3.1]{Po}  we 
introduce the function
\begin{equation*}
\begin{split}
f_{\e}(\bx):=&\int_{\R^n}\ee^{-\e^2 |\btt|^2/4}\, \sigma(\btt) 
\, \ee^{-i(\btt,\bx)}\, d \btt 
=\sum_{\bj} \sum_{|\bb|=0}^{L_\bj}v_{\bj,\bb}\,\Ss_{\bb}(\partial_\bx)
\int_{\R^n}\ee^{- \e^2 |\btt|^2/4} 
\, \ee^{\, i(\btt,\bx_\bj- \bx)}\, d \btt \\
=& \, \e^{-n} \sum_{\bj}
\sum_{|\bb|=0}^{L_\bj}v_{\bj,\bb}\,\Ss_{\bb}(\partial_\bx)
\int_{\R^n}\ee^{- |\btt|^2/4} 
\, \ee^{\,i(\btt,\bx_\bj- \bx)/\e}\, d \btt \\
=&\Big(\frac{4 \pi}{\e^2}\Big)^{n/2} 
\sum_{\bj}
\sum_{|\bb|=0}^{L_\bj}v_{\bj,\bb} \, \Ss_{\bb}(\partial_\bx)
\ee^{- |\bx- \bx_\bj|^2/ \e^2} \, . 
\end{split}
\end{equation*} 
Let us fix $\bk$ and  
consider $f_{\e}(\bx)$ for $|\bx - \bx_\bk| < \e$   
for sufficiently small $\e> 0$.
We have $\Ss_{\bb}(\partial_\bx)\ee^{- |\bx - \bx_\bj|^2/ \e^2} \to 0$
as $\e \to 0$
and
\[
\sum_{|\bb|=0}^{L_\bk} v_{\bk,\bb} \, \Ss_{\bb}(\partial_\bx)
\ee^{- |\bx- \bx_\bk|^2/ \e^2} 
= \sum_{|\bb|=0}^{L_\bk}  v_{\bk,\bb}\Ss_{\bb}( \e^{-1} \partial_\btt)
\ee^{- |\btt|^2} \Big|_{\btt=(\bx- \bx_\bk)/\e} \, .
\]
Because of 
$f_{\e}(\bx) = 0$ for all $\e>0$
there exist $\e_0$
such that
\begin{equation} \label{pzero}
\sum_{|\bb|=0}^{L_\bk} v_{\bk,\bb} \, \Ss_{\bb}(\partial_\bx)
\ee^{- |\bx- \bx_\bk|^2/ \e^2} = 0
\end{equation}
for $\e \le \e_0$. On the other hand,
\[
\sum_{|\bb|=0}^{L_\bk} v_{\bk,\bb} \, \Ss_{\bb}(\partial_\bx)
\ee^{- |\bx- \bx_\bk|^2/ \e^2} =
\e^{-2L_\bk} \ee^{- |\bx- \bx_\bk|^2/ \e^2} \, \Pi_{2L_\bk}(\e) \, ,
\]  
where $\Pi_{2L_\bk}(\e)$ is a polynomial of degree $2L_\bk$ in $\e$ with
coefficients depending on  $|\bx- \bx_\bk|$. Therefore
\eqref{pzero} holds for any $\e> 0$, in particular 
\[
\sum_{|\bb|=0}^{L_\bk} v_{\bk,\bb} \, \Ss_{\bb}(\partial_\bx)
\ee^{- |\bx- \bx_\bk|^2} = 0 \, .
\]
Since by \eqref{delta}
\[
\Ss_{\bb}(\partial_\bx) \ee^{-|\bx- \bx_\bk|^{2}}
= (\bx- \bx_\bk)^{\bb} \ee^{-|\bx- \bx_\bk|^{2}} \, ,
\]
we conclude $v_{\bk,\bb}=0$ for all $\bb$. \hfill $\blacksquare$
\bigskip

Let now for given $\Sigma$ and degrees $L_{\bj}$
the coefficient vector $\cc=\{c_{\bj,\bb}\}$
be a unique solution of the linear system \eqref{sys1}.
To estimate $r = Q(\cc)$ we denote by
$\by_{\bmu} \in \Sigma$ the point closest to $\boldsymbol{0}$
and by $L_{\bmu}$  the degree of the polynomial $\Pol_{\bmu}$.
\begin{lemma}\label{teoremesti}
The minimal value of \eqref{defrm}
can be estimated by
\begin{equation*}
r \le  
\Big(\frac{\pi}{2}\Big)^{n/2} \frac{D^{L_{\bmu}+1+n/2}\, 
|\by_{\bmu}|^{2(L_{\bmu}+1)}}
{D_0^{2(L_{\bmu}+1)}(L_{\bmu}+1)!} \, .
\end{equation*}
\end{lemma}
\textit{Proof.} 
It follows from the representation \eqref{ast1n} that
\begin{equation*}
\begin{split}
r&=\int_{\R^n} 
 \Big(\sum_{\by_\bj \in \Sigma} \ee^{\,(D-D_0)|\by_{\bj}|^2/D_0^2}  
\sum_{|\bb|=0}^{L_{\bj}} 
   c_{\bj,\bb} T_{\bb}(\btt-\by_\bj)
  \,\ee^{-|\btt-\btt_{\bj}|^2/D}
 - \ee^{-|\btt|^2/D} \Big)^2 \, d \btt \\
&\le \min_{\Pol \in \Pi_{L_{\bmu}}} \int_{\R^n} 
\Big( \Pol(\btt) \ee^{-|\btt-\btt_\bmu|^2/D}
- \ee^{-|\btt|^2/D} \Big)^2 \, d \btt \\
&= \Big(\frac{D}{2}\Big)^{n/2}\min_{{\Pol \in \Pi_{L_{\bmu}}}} 
\int_{\R^n} \ee^{-|\btt|^2} \Big(\Pol(\btt) -
\ee^{ - |\bz_{\bmu}|^2}\ee^{- \sqrt{2} (\btt,\bz_{\bmu})}
\Big)^2 \, d \btt
 \end{split}
 \end{equation*}
with $\btt_\bmu= D \by_{\bmu}/D_0$,
$\bz_\bmu= \sqrt{D} \by_{\bmu} /D_0$,
and $\Pi_{L_{\bmu}}$ denotes the set of polynomials
of degree $L_{\bmu}$.
The minimum
is attained when
\[
\Pol(\btt)=\frac{1}{\sqrt{2^{|\bb|} \bb! \pi^{n/2}}}
\sum_{|\bb|=0}^{L_{\bmu}}a_{\bb} H_\bb(\btt)
\]
with the coefficients
\begin{equation*}
\begin{split}
a_{\bb}&=\frac{\ee^{-|\bz_{\bmu}|^2}}{\sqrt{2^{|\bb|}\, \bb! \,\pi^{n/2}}}
\int_{\R^n}{\ee}^{-|\btt|^2}H_\bb(\btt)
\ee^{- \sqrt{2}(\btt,\bz_{\bmu})}
\, d\btt
=\frac{{\ee}^{-|\bz_{\bmu}|^2}}{\sqrt{2^{|\bb|}\bb! \pi^{n/2}}}
\int_{\R^n}{\ee}^{-\sqrt{2}(\btt,\bz_{\bmu})} 
( -\partial_\btt)^{\bb}
{\ee}^{-|\btt|^2} \> d \btt \, .
\end{split}
\end{equation*}
Integrating by parts, we obtain
\[
a_{\bb}=\pi^{n/4}
\frac{(-1)^{|\bb|}\bz_{\bmu}^\bb}{\sqrt{\bb!}}\,
{\ee^{-|\bz_{\bmu}|^2/2}} \, ,
\]
which together with
\[
\sum_{|\bb|=L_{\bmu}+1}^{\infty}\frac{\bz_{\bmu}^{2
\bb}}{\bb!}=\sum_{s=L_{\bmu}+1}^\infty
\frac{|\bz_{\bmu}|^{2s}}{s!}\leq
\frac{|\bz_{\bmu}|^{2(L_{\bmu}+1)}}{(L_{\bmu}+1)!} \, 
{\ee}^{|\bz_{\bmu}|^2}
\]
leads to
\begin{equation*}
\qquad r \le \Big(\frac{D}{2}\Big)^{n/2} 
\sum_{|\bb|=L_{\bmu}+1}^{\infty} 
|a_{\bb}|^2 =
\pi^{n/2} \Big(\frac{D}{2}\Big)^{n/2} 
\frac{|\bz_{\bmu}|^{2(L_{\bmu}+1)}}{(L_{\bmu}+1)!} \, . 
\qquad \blacksquare
\end{equation*}

\subparagrafo{Approximate partition of unity with Gaussians}
Now we are in position to prove the main result of this section.
Suppose that the nodes $\{ \bx_\bj \}$ 
are as described in subsection~\ref{SS:61} and let 
$G_1 \cup G_2$ be the
associated piecewise uniform grid with stepsizes $h_1$ and $h_2$.
Assign to each grid point $g_\ell \in G_\ell$, $\ell = 1,2$,
a finite set of nodes $\Sigma(g_\ell)$,
fix a common degree $L$ for all polynomials $\Pol_\bj$ in  \eqref{th0}
and solve the linear system
\begin{equation}\label{sys2}
\sum_{\bx_\bj \in \Sigma(g_\ell)} \sum_{|\bb|=0}^{L} 
\C_{\bb,\bg}\Big(\frac{\bx_{\bj}-g_\ell}{h_\ell} ,
\frac{\bx_{\bk}-g_\ell}{h_\ell}\Big) c_{\bj,\bb}(g_\ell)=
\C_{\boldsymbol{0},\bg}\Big(\boldsymbol{0},
\frac{\bx_{\bk}-g_\ell}{h_\ell} \Big) 
\end{equation}
for all $\bx_{\bk} \in \Sigma(g_\ell)$ and 
$0 \le |\bg| \le L$ with 
\[
\C_{\bb,\bg}(\bx,\by) = 
\Ss_{\bb}(- \dd \partial_{\bx})
\Ss_{\bg}(-\dd \partial_{\by}) \,
\ee^{\,(D-D_0)(|\bx|^2+|\by|^2)/D^2_0} \, 
\ee^{ -D|\bx-\by|^2/2 D^2_0} \, ,
\]
$D_0< D$ is some arbitrary positive number.
Following \eqref{lsdef} define the polynomials 
\begin{equation} \label{lsdef1}
\begin{split}
&\Pol_{\bj}
\Big( \frac{\bx-\bx_\bj}{h_1 \dd} \Big)
=\sum_{g_1 \in G(\bx_\bj)}
\sum_{|\bb|=0}^{L}
c_{\bj,\bb}(g_1) \Big( \frac{\bx-\bx_\bj}{h_1\dd} \Big)^\bb
\, , \qquad  \quad \bx_\bj \in J_1  , \\
&\Pol_{\bk}
\Big( \frac{\bx-\bx_\bk}{h_2 \dd} \Big)
=\sum_{g_2 \in G(\bx_\bk)}
\sum_{|\bb|=0}^{L} \tilde a_{g_2} 
c_{\bk,\bb}(g_2) \Big( \frac{\bx-\bx_\bk}{h_2\dd} \Big)^\bb
\, , \hskip5.3mm \bx_\bk \in J_2  , 
\end{split}
\end{equation}

\begin{thm}
Under Conditions~{\em \ref{cond1}} and {\em \ref{cond4}} 
on the scattered nodes $\{\bx_\bj\}$ for any $\e>0$ there exist $D>0$ and $L$ 
such that the function \eqref{th0} 
is an approximate partition of unity satisfying  
\[
|\Th(\bx)-1|<\e \quad \mbox{ for all } \,  \bx \in \R^n \,,
\]
if the polynomials $\{\Pol_{\bj}\}$ of degree $L$
are generated via \eqref{lsdef1}
by the solutions
$\{c_{\bj,\bb}(g_\ell)\}$ of the linear systems \eqref{sys2}
for all $g_\ell \in G_1 \cup G_2$. 
\end{thm}
\textit{Proof.} 
From \eqref{thm1} we have to show
that 
\begin{equation} \label{thm110}
\sup_{\R^n} \Big(
\sum_{g_1 \in G_1}
\Big|\o_{g_1} \big(\frac{\bx}{h_1}\big)\Big| +
\sum_{g_2 \in G_2} \tilde a_{g_2} 
\Big|\o_{g_2} \big(\frac{\bx}{h_2}\big)\Big| \Big)
\le \frac{\e}{2}
\end{equation}
if $L$ is sufficiently large.
We start with estimating the first sum
\[
\sum_{g_1 \in G_1}
\Big|\o_{g_1} \big(\frac{\bx}{h_1}\big)\Big| \, ,
\]
where $g_1= h_1 \bm$, $\bm \in Z_1 \subset \Z^n$.
Using \eqref{omegam_2} we can write
\[
\o_{g_1} \big(\frac{\bx}{h_1}\big) =
\o \big(\frac{\bx}{h_1}-\bm\big) \, ,
\]
where the points $\by_\bj$ in \eqref{omegam_2}
are given by $\by_\bj = \bx_\bj/h_1-\bm$, $\bx_\bj \in \Sigma(g_1)$.
By Lemmas \ref{constr} and \ref{teoremesti} we have
\begin{equation*}
\sum_{g_1 \in G_1}
\Big|\o_{g_1} \big(\frac{\bx}{h_1}\big)\Big| 
\le c_3 \Big(\frac{\pi}{2}\Big)^{n/4}
\sum_{\bm \in Z_1} \ee^{- \rho|\bx/h_1-\bm|^2} 
 \frac{D^{(L_{\bmu_\bm}+1+n/2)/2}} 
{D_0^{L_{\bmu_\bm}+1}\sqrt{(L_{\bmu_\bm}+1)!}} \, 
\Big|\frac{\bx_{\bmu_\bm}}{h_1} - \bm\Big|^{L_{\bmu_\bm}+1} 
\end{equation*}
where 
$\bx_{\bmu_\bm}\in \Sigma(g_1)$ is the node
closest to $g_1= h_1 \bm$
and $L_{\bmu_\bm}$ is the degree of the polynomial
$\Pol_{\bmu_\bm,g_1}$. 
Since $|\bx_{\bmu_\bm} - h_1 \bm| \le \kappa_1 h_1$ by Condition~\ref{cond1}
and $L_{\bmu_\bm}=L$ for all $\bmu_\bm$
we conclude that
\begin{equation} \label{sumest}
\sum_{g_1 \in G_1}
\Big|\o_{g_1} \big(\frac{\bx}{h_1}\big)\Big| 
\le c_3 \Big(\frac{\pi}{2}\Big)^{n/4}
 \frac{D^{(L+1+n/2)/2}\, \kappa_1^{L+1} } 
{D_0^{L+1}\sqrt{(L+1)!}} 
\, \sup_{\R^n} \sum_{\bm \in Z_1} \ee^{- \rho|\bx/h_1-\bm|^2} \, .
\end{equation}
From
\[
  \rho = \frac{D-D_0}{(D-D_0)^2+DD_0} 
\in (0,D) \, \quad \mbox{for any fixed}  \quad  D_0 \in (0,D) \, ,
\]
we see, that for fixed $D$ and $D_0$
\begin{equation} \label{firtpart}
\sum_{g_1 \in G_1}
\Big|\o_{g_1} \big(\frac{\bx}{h}\big)\Big| \to 0 \quad \mbox{if }
\quad  L \to \infty \, .
\end{equation} 

We turn  to 
\[
\sum_{g_2 \in G_2} \tilde a_{g_2} 
\Big|\o_{g_2} \big(\frac{\bx}{h_2}\big)\Big| 
\]
with $g_2=h_1\bm + h_2 \bk, \,\bm \in Z_2, \bk \in S$.
Using \eqref{omegam_2} we have
\[
\o_{g_2} \big(\frac{\bx}{h_2}\big) =
\o \big(\frac{\bx- \bm h_1}{h_2}-\bk\big) \, ,
\]
and the points $\by_\bj$ in \eqref{omegam_2}
are given by $\by_\bj = (\bx_\bj-\bm h_1)/h_2-\bk$ with
$\bx_\bj \in \Sigma(g_2)$.
Hence
\begin{align*}
&\sum_{g_2 \in G_2} \tilde a_{g_2} 
\Big|\o_{g_2} \big(\frac{\bx}{h_2}\big)\Big| 
= \sum_{\bm \in Z_2} \sum_{\bk \in S}
a_\bk \Big| \o \big(\frac{\bx- \bm h_1}{h_2}-\bk\big) \Big|\\
&\le c_3 \Big(\frac{\pi}{2}\Big)^{n/4}
\sum_{\bm \in Z_2} \sum_{\bk \in S}
a_\bk \ee^{- \rho|(\bx-\bm h_1)/h_2-\bk|^2} 
 \frac{D^{(L_{\bmu_\bk}+1+n/2)/2}} 
{D_0^{L_{\bmu_\bk}+1}\sqrt{(L_{\bmu_\bk}+1)!}} 
\Big|\frac{\bx_{\bmu_\bk}-\bm h_1}{h_2} - \bk\Big|^{L_{\bmu_\bk}+1} \, .
\end{align*}
Here $\bx_{\bmu_\bk}\in \Sigma(g_2)$ is the node
closest to $g_2= h_1\bm + h_2 \bk$
and $L_{\bmu_\bk}$ is the degree of the polynomial
$\Pol_{\bmu_\bk,g_2}$. 
By Condition~\ref{cond4} for fixed $D$ and $D_0$
\[
 \frac{D^{(L+1+n/2)/2}} 
{D_0^{L+1}\sqrt{(L+1)!}} 
\Big|\frac{\bx_{\bmu_\bk}-\bm h_1}{h_2} - \bk\Big|^{L+1} 
\le \delta(L) \to 0 \quad \mbox{if }
\quad  L \to \infty 
\]
uniformly for all $g_2 \in G_2$.
Hence we obtain 
\begin{equation} \label{secondpart}
\sum_{g_2 \in G_2} \tilde a_{g_2} 
\Big|\o_{g_2} \big(\frac{\bx}{h_2}\big)\Big|
\le  C_1 \delta(L) \sum_{\bm \in Z_2} \sum_{\bk \in S}
a_\bk \ee^{- \rho|(\bx-\bm h_1)/h_2-\bk|^2} 
\end{equation} 
because of $L_{\bmu_\bk}=L$ for all $\bmu_\bk$. The sum 
\[
 \sum_{\bk \in \Z^n }a_\bk \ee^{- \rho|(\bx-\bm h_1)/h_2-\bk|^2}
=  \Big(\frac{h_1^2} { \pi D(h_1^2- h_2^2)} \Big)^{n/2}
\sum_{\bk \in \Z^n } \ee^{- h_2^2|\bk|^2/ (h_1^2-h_2^2)D}
\, \ee^{- \rho|\bx-\bm h_1 - h_2\bk|^2/h_2^2}
\]
can be easily estimated by using equation \eqref{refin}.
Setting
\[
(h_1^2-h_2^2)D = h_1^2 D_1-h_2^2/\rho
\]  
we derive
\[
D_1= D + \frac{h_2^2}{h_1^2}\Big( \frac{1}{\rho} - D \Big)= D+
H^2 \,  \frac{D_0^2 }{D-D_0} \, .
\]
and after some algebra
\begin{align*}
\sum_{\bk \in \Z^n }& \ee^{- h_2^2|\bk|^2/ (h_1^2-h_2^2)D}
\, \ee^{- \rho|\bx - h_2\bk|^2/h_2^2} \\
&= \Big(\frac{\pi D(1-H^2)}{\rho D_1}\Big)^{n/2}
\ee^{-|\bx|^2/h_1^2 D_1}\Big(1 + O(\ee^{- \pi^2 D^2(1-H^2)/D_1})\Big).
\end{align*}
Therefore we obtain
\[
\sup_{\R^n} \sum_{\bm \in Z_2} \sum_{\bk \in S}a_\bk \ee^{- \rho|(\bx-\bm h_1)/h_2-\bk|^2}
\le C_2 \sup_{\R^n}  \sum_{\bm \in Z_2} \ee^{-|\bx-\bm h_1|^2/h_1^2 D_1}
\le C_3
\]
with some constant $C_3$ depending on $D$, $D_0$ and the space
dimension $n$. 
Now \eqref{thm110} follows immediately from 
\eqref{firtpart} and \eqref{secondpart}.
\hfill  $\blacksquare$

\begin{rem}
{ \em
It can be seen from  \eqref{sumest} that
in principle the parameter $D_0$ can 
be any value of the interval $(0,D)$ not too close the its end
points. In numerical experiments we have not seen 
any significant dependence on this parameter. 
The choice $D_0 = D/2$ might be advantageous because
the differential expressions $\C_{\bb,\bg}(\bx,\by)$ 
simplify to
\[
\C_{\bb,\bg}(\bx,\by) = 
\Ss_{\bb}(- \dd \partial_{\bx})
\Ss_{\bg}(-\dd \partial_{\by}) \,
\ee^{4 (\bx,\by)/D} \, .
\]
}
\end{rem}

\begin{figure}[p] \begin{center}
 \epsfig{file=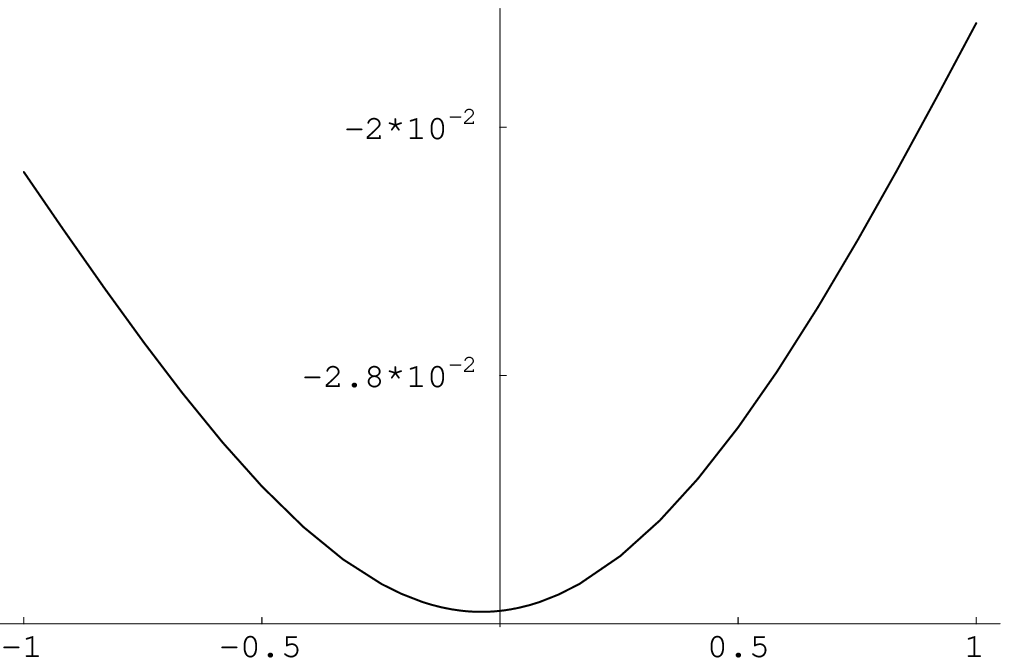,width=3in,height=1.8in}
 \epsfig{file=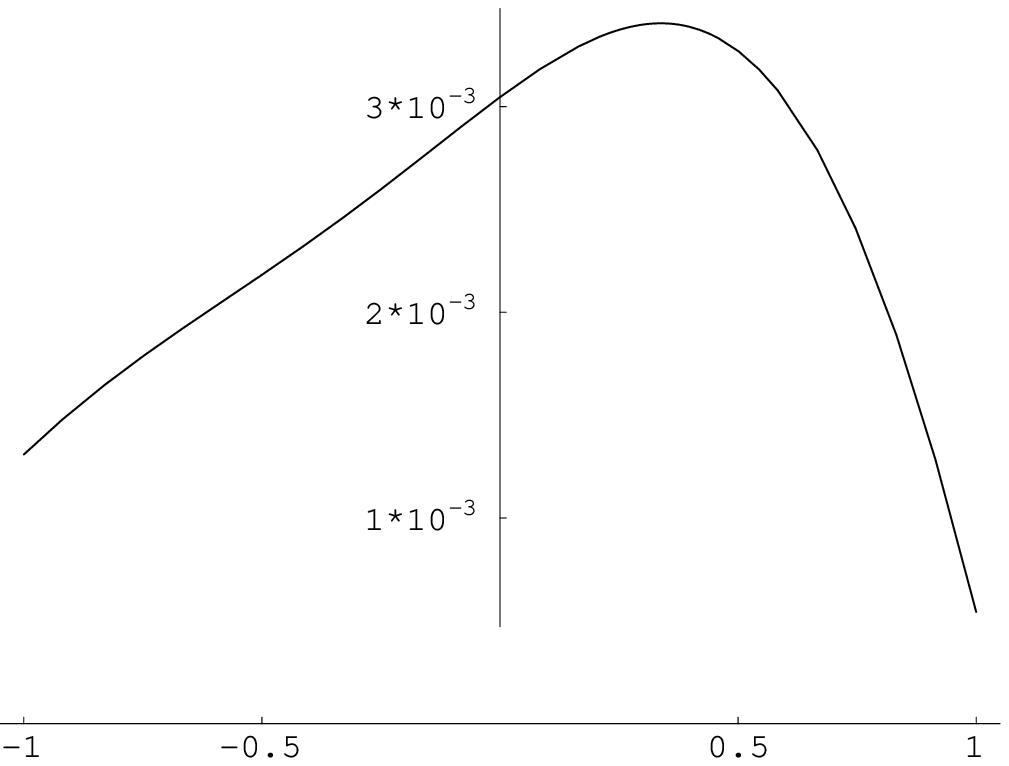,width=3in,height=1.8in}
 \caption{{\small The graph  of $\Th(\bx)-1$ when 
$\Sigma(m)$ consists of $1$ point,  $L=1$ (on the left) and $L=2$ 
(on the right).}}
 \label{fig(1)}
 \epsfig{file=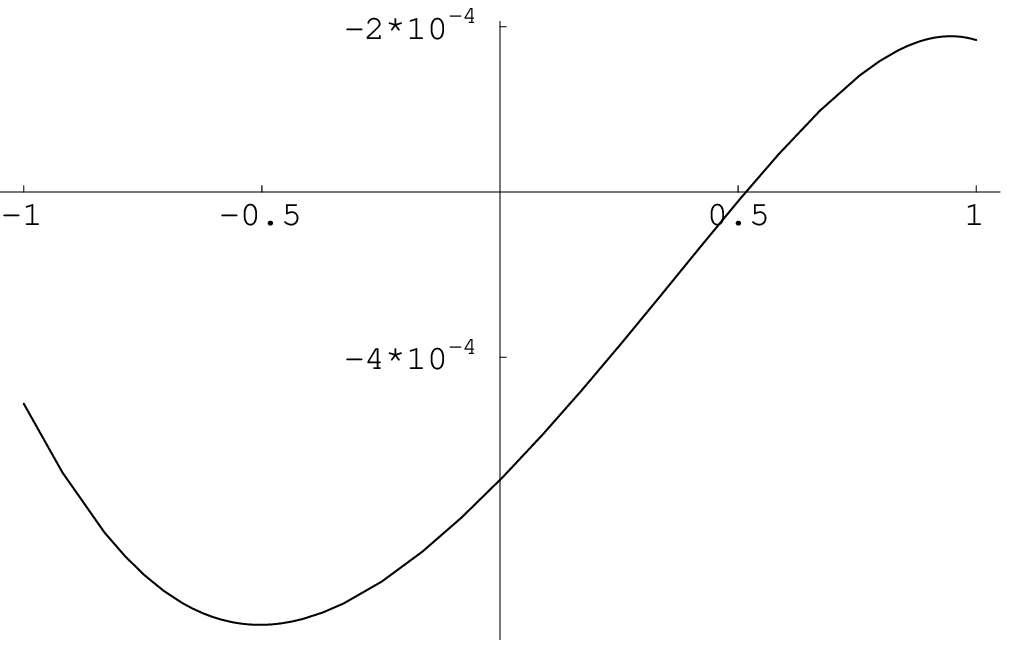,width=3in,height=1.8in}
 \epsfig{file=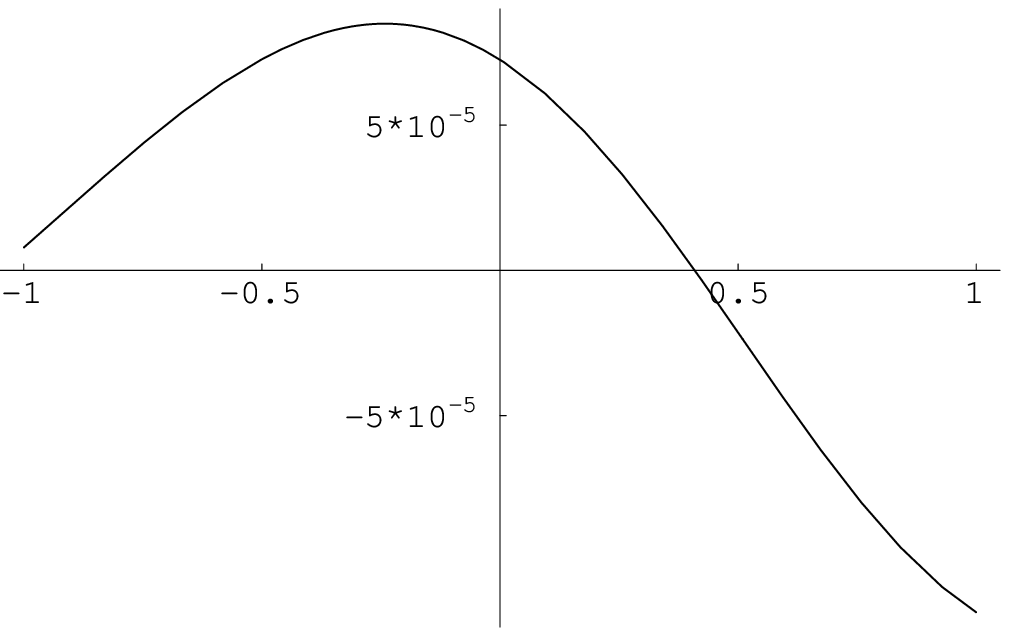,width=3in,height=1.8in}
 \caption{{\small The graph  of $\Th(\bx)-1$ when 
$\Sigma(m)$ consists of $1$ point,  $L=3$ (on the left) and $L=4$ 
(on the right).}}
 \label{fig(2)}
 \epsfig{file=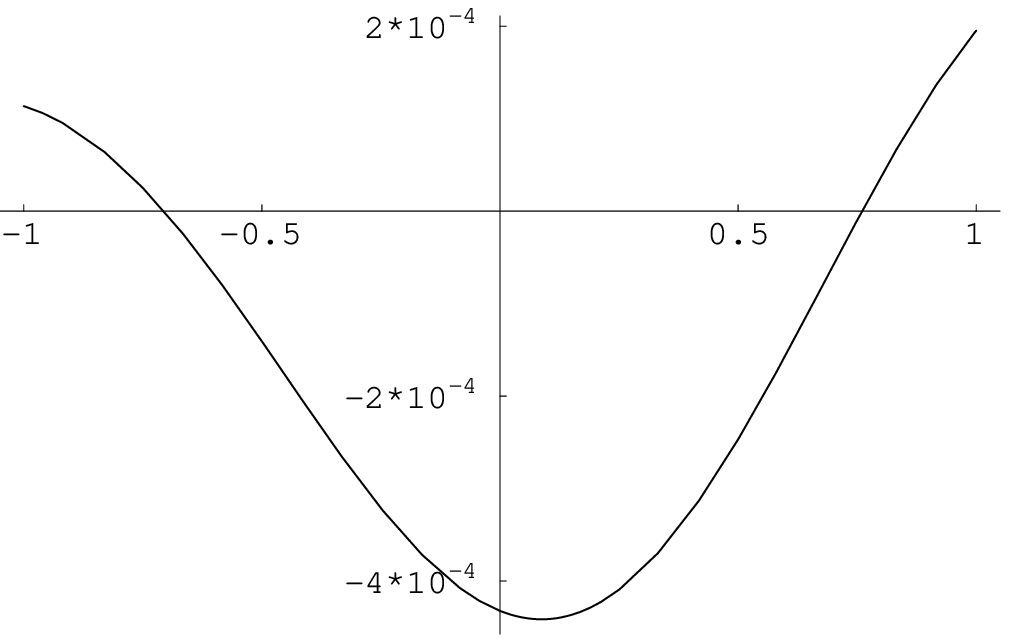,width=3in,height=1.8in}
 \epsfig{file=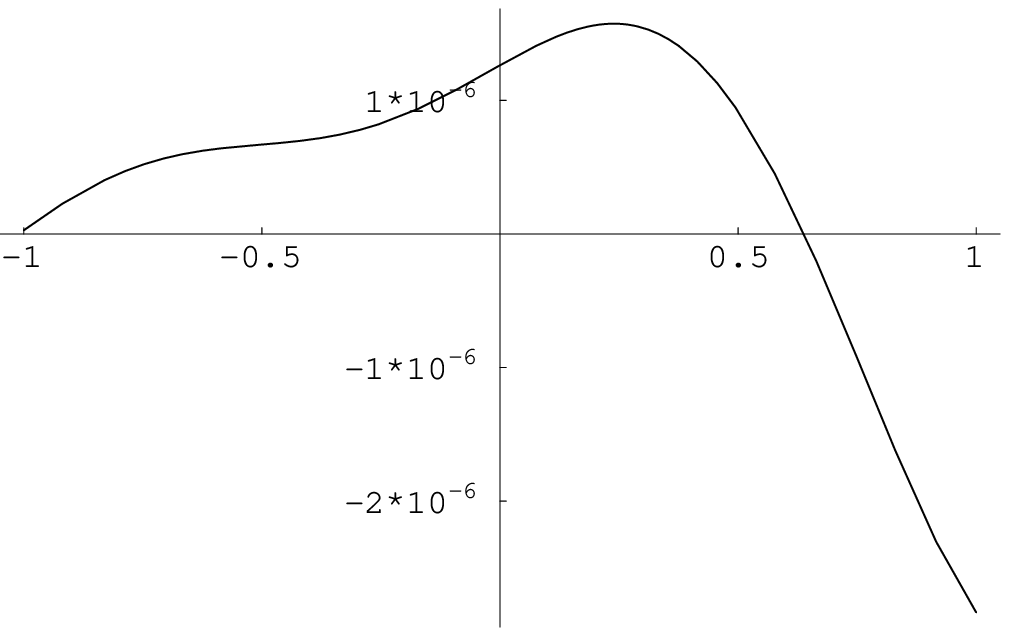,width=3in,height=1.8in}
 \caption{{\small The graph  of $\Th(\bx)-1$ when 
$\Sigma(m)$ consists of $3$ points,  $L=1$ (on the left) and $L=2$ 
(on the right).}}
 \label{fig(4a)}
 \epsfig{file=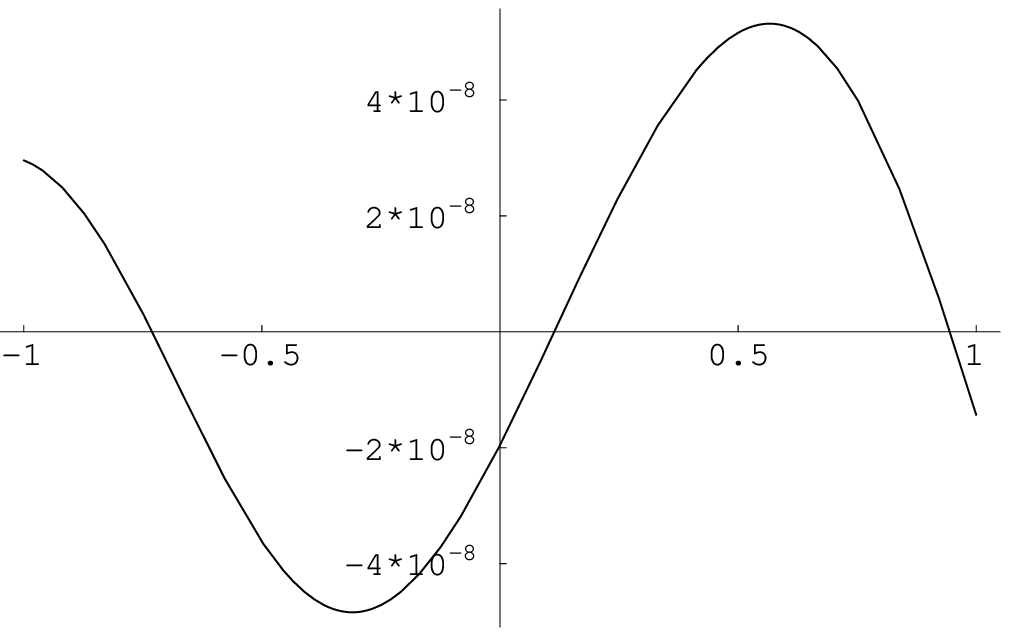,width=3in,height=1.8in}
  \epsfig{file=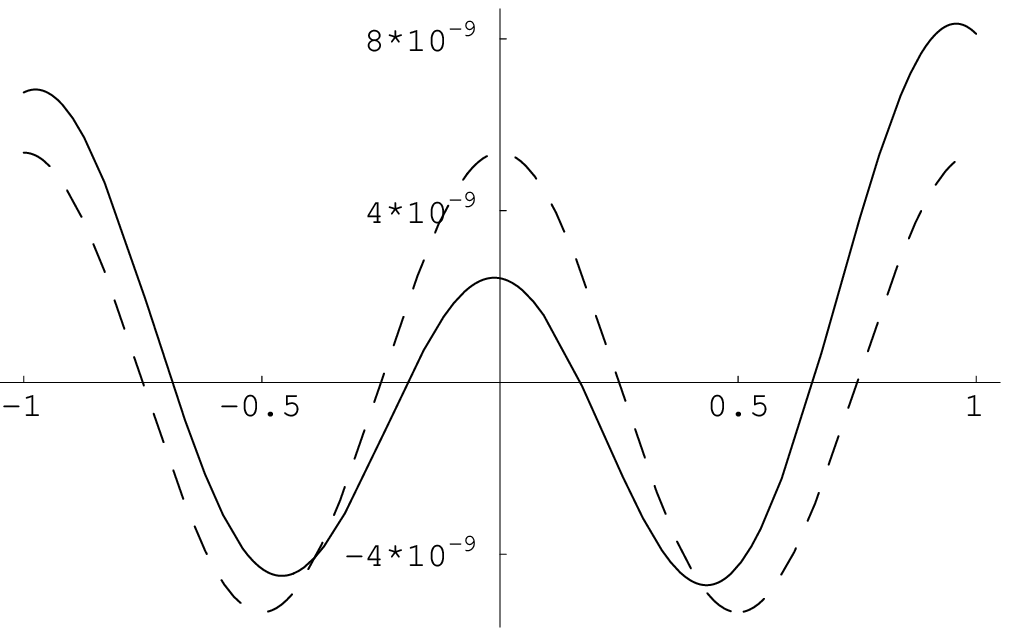,width=3in,height=1.8in}
 \caption{{\small The graph  of $\Th(\bx)-1$ when 
$\Sigma(m)$ consists of $3$ points, $L=3$ (on the left) and $L=4$
(on the right). The saturation term \eqref{defsat} is depicted by dashed lines.}}
 \label{fig(4c)}
 \end{center}
 \end{figure}

\begin{figure}[h] \begin{center}
 \epsfig{file=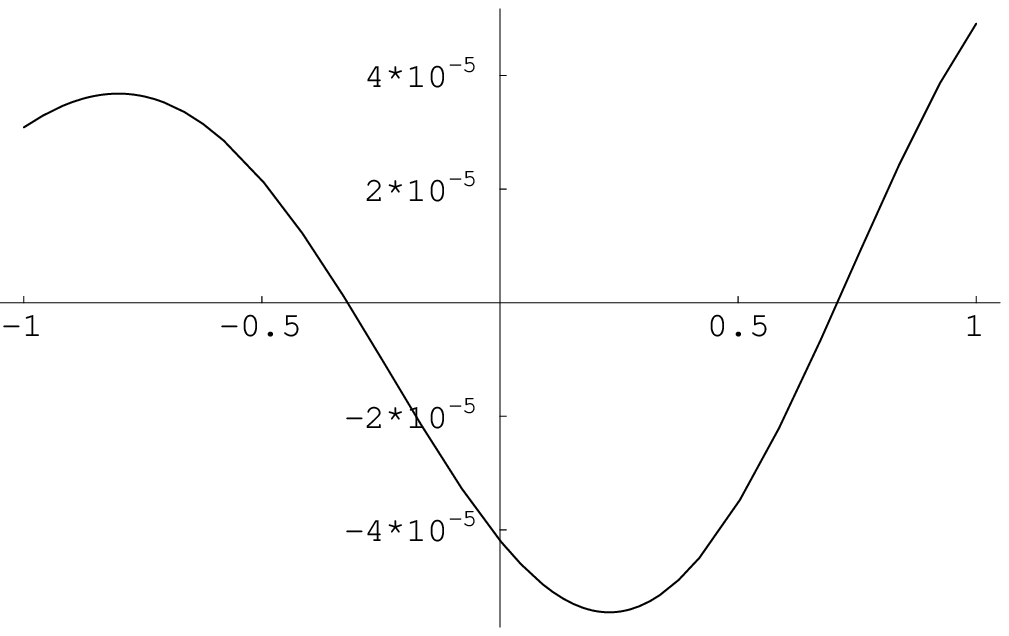,width=3in,height=1.8in}
  \epsfig{file=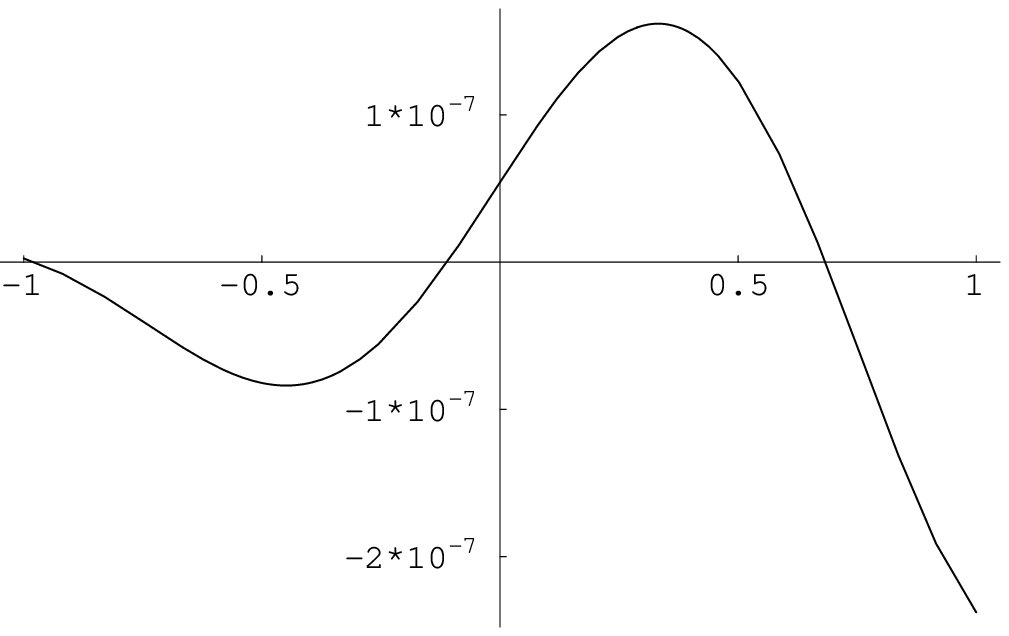,width=3in,height=1.8in}
 \caption{{\small The graph  of $\Th(\bx)-1$ when 
$\Sigma(m)$ consists of $5$ points, $L=1$ (on the left) and $L=2$
(on the right).}}
 \label{fig(4e)}
 \epsfig{file=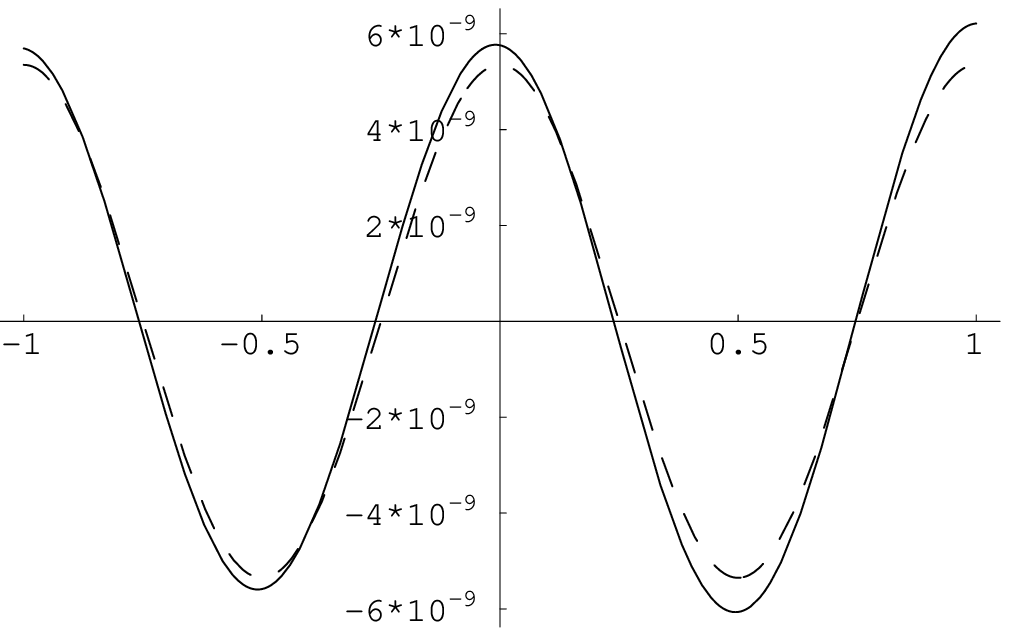,width=3in,height=1.8in}
  \epsfig{file=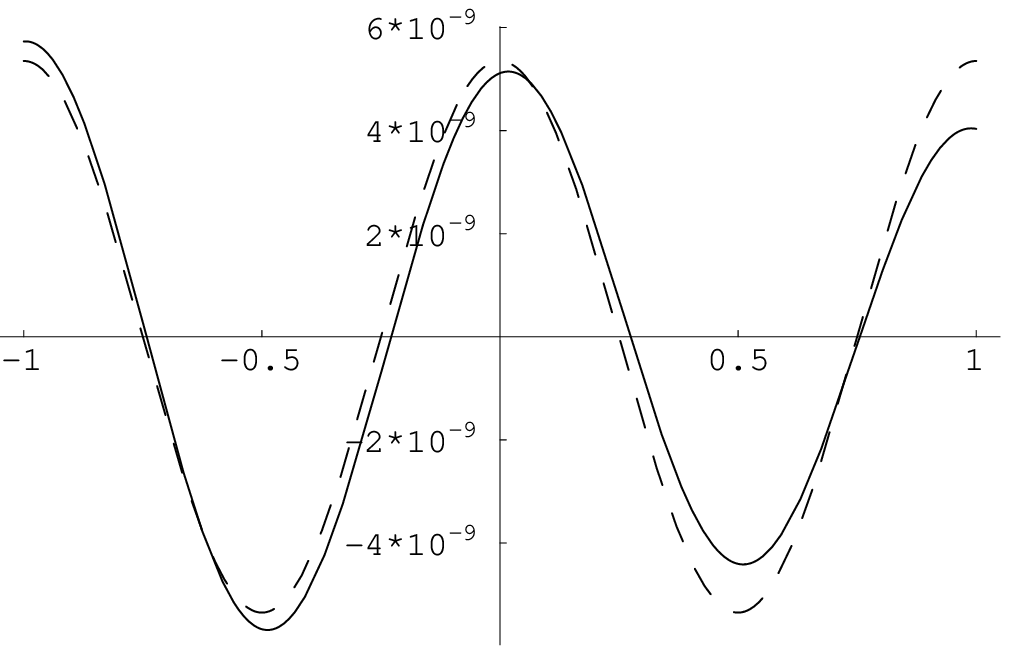,width=3in,height=1.8in}
 \caption{{\small The graph  of $\Th(\bx)-1$ when 
$\Sigma(m)$ consists of $5$ points, $L=3$ (on the left) and $L=4$
(on the right). The saturation term \eqref{defsat} is depicted by dashed lines. 
}}
 \label{fig(4g)}
 \end{center}
 \end{figure}

\subparagrafo{Numerical Experiments}\label{S:NUM2}
We have tested the construction (\ref{lsdef1}, \ref{sys2})
for a quasi-uniform distribution of nodes on $\R$
with the parameters $D=2$, $D_{0}=3/2$, $h=1$, $\k_{1}=1/2$. 
To see the dependence of the approximation error
from the number of nodes in $\Sigma(m)$, $m \in \Z$, and the degree of polynomials
we provide graphs of the difference to $1$ for the following cases :
\begin{itemize}
\item[-] $\Sigma(m)$ consists of $1$ point, $L=1,2$ (Fig.
\ref{fig(1)})  and $L=3,4$ (Fig.
\ref{fig(2)});
\item[-] $\Sigma(m)$ consists of $3$ points, $L=1,2$ (Fig.
\ref{fig(4a)}) and $L=3,4$ (Fig.
\ref{fig(4c)});
\item[-] $\Sigma(m)$ consists of $5$ points, $L=1,2$ (Fig.
\ref{fig(4e)}) and $L=3,4$ (Fig.
\ref{fig(4g)}).
\end{itemize}
As expected, the approximation becomes better
with increasing degree $L$ and  
more points in the subsets $\Sigma(m)$.
The use of only one node in $\Sigma(m)$
reduces the  approximation error by a factor $10^{-1}$
if $L$ increases by $1$.
The cases of $3$ and $5$ points indicate, that 
enlarging the degree $L$ of the polynomials by $1$  gives
a factor $10^{-2}$ for the approximation error.
 
One should notice, 
that the plotted total error 
consists of two parts. Using (\ref{lsdef1}, \ref{sys2})
we approximate the $\theta$-function
    \begin{equation} \label{defpos} 
(2 \pi)^{-1/2}\sum_{m \in \Z} \ee^{-(x-m)^2/2}= 1+
2 \sum_{j=1}^{\infty} \ee^{-2 \pi^2 j^2} \cos 2 \pi j x \, .
    \end{equation}
Hence, the plotted total error 
is the sum of the difference between \eqref{th0} and \eqref{defpos}
and the function
    \begin{equation} \label{defsat} 
2 \sum_{j=1}^{\infty} \ee^{-2 \pi^2 j^2} \cos 2 \pi j x \, ,
    \end{equation}
which is the saturation term obtained on the uniform grid.
The  error plots in Figure \ref{fig(4c)}
for $L=4$ and in Figure \ref{fig(4g)}
show that the total error is already majorized by \eqref{defsat}, 
which is shown by dashed lines.
\vskip15mm


\end{document}